\documentclass[EJP]{ejr} 

\usepackage{graphicx}
\usepackage{nccmath}
\usepackage[section]{placeins}
\usepackage[para]{threeparttable}
\usepackage{amsthm}
\theoremstyle{plain} 

\SHORTTITLE{High, Low and Close in Brownian Estimation} 

\TITLE{The Value of the High, Low and Close in the Estimation of Brownian Motion: Extended Version} 

\AUTHORS{%
  Kurt Riedel\footnote{Millennium Partners LLC
    \EMAIL{kurt.riedel@gmail.com}}}

\date{\today}

\AMSSUBJ{60J65, 60J70, 91G60} 

\SUBMITTED{August 14, 2019, Revised Feb.~12, 2020, Oct.~27, 2020.
This extended version differs from the SISP article by including Section \ref{TDCHL}
and additional plots in Section \ref{SCVsect}. } 

\ARXIVID{1911.05280} 

\VOLUME{0}
\YEAR{2020}
\PAPERNUM{0}
\DOI{10.1214/YY-TN}

\ABSTRACT{
The conditional density of Brownian motion is considered given the max, $B(t|\max)$, as well as those with additional information:
  $B(t|close, \max)$, $B(t|close, \max, \min)$ where the close is the final value:
  $B(t=1)=c$ and $t \in [0,1]$.
  The conditional expectation and conditional variance of Brownian motion are evaluated subject
  to one or more of the statistics: the close (final value), the high (maximum),  the low (minimum). 
  Computational results displaying both the expectation and variance in time are presented
  and compared with the theoretical values. We tabulate the time averaged
  variance of Brownian motion conditional on knowing various extremal properties of the motion.
  The final table shows that knowing the high is more useful than knowing the final value among other results.
  Knowing the open, high, low and close reduces the time averaged variance to $42\%$ of the value of knowing only the open and close (Brownian bridge). 
} 

\newcommand{\FSCL}{.5}

\newcommand{\close}{\mathrm{close}}

\newcommand{\BAL}{\begin{align}}
\newcommand{\NAL}{\end{align}}
\newcommand{\BEQ}{\begin{equation}}
\newcommand{\EEQ}{\end{equation}}
\newcommand{\NEQ}{\end{equation}}
\newcommand{\Beq}{\begin{equation*}}
\newcommand{\Neq}{\end{equation*}}
\newtheorem{lemma}[theorem]{Lemma} 

\newcommand{\sig}{\sigma}
\newcommand{\sgg}{\sigma^2}
\newcommand{\sggh}{\hat{\sigma}^2}
\newcommand{\sggt}{\sigma^2_t}

\newcommand{\jh}{\hat{j}}
\newcommand{\kh}{\hat{k}}
\newcommand{\qb}{\bar{q}}
\newcommand{\Ab}{\bar{A}}
\newcommand{\Bb}{\bar{B}}

\newcommand{\Btb}{\bar{\tilde{B}}}

\newcommand{\Xh}{\hat{X}}
\newcommand{\gtl}{\tilde{g}}
\newcommand{\Gam}{\Gamma}
\newcommand{\Gamt}{\tilde{\Gamma}}
\newcommand{\ah}{\hat{a}}

\newcommand{\phisg}{\phi_{\sigma^2}}
\newcommand{\phitsg}{\phi_{t\sigma^2}}
\newcommand{\phihsg}{\phi_{\hat{t}\sigma^2}}
\newcommand{\sgp}{\sigma^2}

\newcommand{\vt}{\tilde{v}}
\newcommand{\wt}{\tilde{w}}
\newcommand{\Bc}{\mathcal{B}}
\newcommand{\taut}{\tilde{\tau}}
\newcommand{\tauh}{\hat{\tau}}
\newcommand{\tauth}{\hat{\tilde{\tau}}}
\newcommand{\Prf}{{\em Proof:}\, }

\begin{document}
\date{\today}
\section{Introduction}
\graphicspath{{/home/kriedel/Tex/brownFig/}}

In today's financial markets, every tick is archived. In analyzing events in the ancient past (1970s) or less automated markets like
credit default swaps or emerging market bonds (roughly pre-2013), the only data that typically is available is the open, high, low, close data.
An entire field, chartist analysis, uses these descriptors as the "sufficient statistics" for prediction. 
This paper defines the probability distribution of $B(t|high,low,close)$ and calculates its expectation. 
Our formulas allow us to interpolate the price signal as  $E[B(t|open,high,low, close)]$ over all time in $[0,1]$ given
any data source that only has open, high, low, close data. We think of the  open, high, low and close as ``statistics'' which we will
use to estimate the mean and variance of the process.

Classically, most of the financial forecasting based on charts uses at most four pieces of
information for each day, the opening price (open), the closing price (close), the maximum price (high)
and the minimum price, (low) \cite{Morris}. We address the issue of how much additional information the high and low
carry beyond that of the open and close. We measure the ``value'' by the reduction
of the variance of the Brownian motion given one or both of the high, $h$, and low, $\ell$ .
The variance of the path of a Brownian motion is $V(t) =t$ which integrates to $\int_0^1 V(s) =1/2$.
For the Brownian bridge pinned to $B(t)= c$, the variance is independent of the terminal value, $c$, and satisfies $V(t)=t(1-t)$.
Integrating variance of the Brownian bridge from zero to one yields an average variance, $\int_0^1 V(s) =1/6$. Thus knowledge of both the open
and the close significantly reduces the variance of the process. Our results allow us to calculate the variance of Brownian motion given the
high, low and close, $V(t|h,\ell,c)$.

There have been a number of studies that use the open, high, low and close to improve the estimate of the volatility (standard deviation)
of the Brownian motion \cite{GarmenKlass,McLeish,Meillijson,RogersSatchell}. In contrast, we assume that the variance is given and standardized to $\sgg$.
In reality, the volatility of financial time series are unknown, bursty, and temporally non-uniform on many time scales.
Given a model of the time dependence of the volatility, one must transform time to an equal volatility time.
For this paper, we neglect this difficult problem and proceed with the studying standardized Brownian motion.

Let $B(t)$ be a Brownian motion on $[0,1]$ and $B_c(t)$ be the Brownian motion restricted
to $B_c(t=1)=c$. We allow an arbitrary const variance, $E[B(s)^2]=\sigma^2$.
Our notation tracks the excellent compendium of results by Devroye \cite{Devroye}. Many of the results summarized in Sections  \ref{BEsect} can be found there.
We consider the distribution of $B(t)$ conditioned on one or more of the statistics: $B(t=1)=c$, $\max_{t \in[0,1]}B(t)=h$ and
$\min_{t \in[0,1]}B(t)=\ell$. We evaluate the conditional density of $B(t|close, \max)$ and $B(t|close, \max, \min)$
using Chapman-Kolmogorov type calculations \cite{FPY}. The conditional densities of $B(t|\max)$ and $B(t|\max, \min)$ are found by
integrating the earlier densities over $c$. Our primary goal is to evaluate the conditional mean and conditional variance
of $B(t)$ in these cases.
For several cases, explicit formulas for the moments are given.
The location of the minimum and the location of the maximum are unknown and not used in our analysis.
A number of other studies \cite{Shepp,PY,Imhof3,RiedelArgMx} examine the distribution of Brownian motion given its maximum
and the location of its maximum (Williams-Denisov decomposition). In a sister publication \cite{RiedelArgMx}, we
compute the conditional density and moments of $B(t|close, \max, argmax)$.
In theory, we could integrate the densities over the location of the maximum. Our experience is that this is analytically intractable.

All of our results are for Brownian motion, $B(t)$, on $[0,1]$ with  $Var[B(t)]=\sgp$. We will often
use the notation $\sgp_t \equiv \sgp *t(1-t)$.
Section \ref{BEsect} reviews results on the density/distribution of extrema of Brownian motion.
Section \ref{HCcalcSect} derives analytic formuli for $E[B(t|c,h)]$ and $Var[B(t|c,h)]$.
Section \ref{HLCcalcSect} derives analytic formuli for $E[B(t|c,h,\ell)]$ and $Var[B(t|c,h,\ell)]$.
Section \ref{NumSect} reviews our numerical simulations.
Section \ref{SCVsect} plots $E[B(t|c)]$ and $Var[B(t|c)]$ as well as $E[B(t|h)]$ and $Var[B(t|h)]$. It then computes Feller's distribution
for the range, $\Delta=h-\ell$, and compares with our simulation results. 
Section \ref{TDCH} plots the $E[B(t|c,h)]$ and $Var[B(t|c,h)]$ for a variety of different values of $(c,h)$.
Section \ref{CompCHSect} compares the analytic formuli in Section \ref{HCcalcSect} with the simulation results in Section \ref{TDCH}.
Section \ref{TDCHL} plots the $E[B(t|c,h,\ell)]$ and $Var[B(t|c,h,\ell)]$ for a variety of different values of $(c,h,\ell)$.
Section \ref{CompCHLSect} compares the analytic formuli in Section \ref{HLCcalcSect} with the simulation results. 
Section \ref{HLSect} derives the distribution $p(t|h,\ell)$ by integrating over the closing value in $p(t,h,\ell,c)$.
Section \ref{Conclude} discusses and summarizes results. Especially important are Table 1 
and Figure \ref{fi:AvgVarGivenBIG} as they demonstrate the variance reduction from using the high and low in the
estimation of $B(t)$.

\section{Distributions of Brownian Extrema \label{BEsect}}

The study of Brownian extrema date back to the founders of the field \cite{Levy}. Our brief discussion follows \cite{Devroye} with additional results taken from \cite{Shepp,DI}, \cite{DIM,ImhofHung,Borodin}.
See also \cite{BertoinPitman,KS,PY,Yor}.
We denote the Gaussian density by $\phi_{\sigma^2}(x) = (2\pi \sigma^2)^{-.5} \exp(-x^2/2\sgp)$.
The density of the high, (maximum of $B(t)$), $h$ is that of the half normal: $p(h)=2 \phisg(h) =\sqrt{\frac{2}{\pi \sgp }} \exp(-h^2/2\sgp)$, $h>0$.
The classic result \cite{WilliamsBAMS,ImhofJAP} derived using the reflection principle is
\begin{theorem} \label{refB}
The joint distribution of the close, $c$, the high, $h$ is
\BEQ
P(h,c)= P(\max\{B(s),\ s \in [0,1]\}\le h, B(1)=c) = \phisg(c) - \phisg(2h-c) \ .
\NEQ
The marginal density satisfies
\BEQ \label{eq:denHC}
p(h,c) \equiv p(h=\max\{B(s),\ s \in [0,1]\}, B(1)=c) = \frac{2(2h-c)}{\sgp} \phisg(2h-c) \ , 
\NEQ
where $ h\ge 0$, $h \ge c$.
\end{theorem}
Here $P(h,c)$ is a distribution in $h$ and a density in $c$. 
The conditional density, $p(c|h)$, is given by
\BEQ \label{eq:denCcondH}
p(c|h) =p(h,c)/p(h) =  \ \frac{(2h-c) \exp(-(2h-c)^2/2\sgp)}{\sgg \exp(-h^2/2\sgp)} \ .
\NEQ
Using \eqref{eq:denCcondH}, we find 
\BEQ \label{eq:EVCH}
E[c|h]= h - \sigma \sqrt{\frac{\pi}{2}} erfc(\frac{h}{\sqrt{2\sgp}}) \exp(\frac{h^2}{2\sgp}) \ , \quad E[c^2|h]= h^2 +2\sgp -4h\sigma \sqrt{\frac{\pi}{2}} erfc(\frac{h}{\sqrt{2\sgp}})  \exp(\frac{h^2}{2\sgp}) \ ,
\NEQ
\BEQ
 Var[c|h]= h^2 +2\sgp -2h\sigma\sqrt{\frac{\pi}{2}} erfc(\frac{h}{\sqrt{2\sgp}})  \exp(\frac{h^2}{2\sgp}) - \frac{\pi\sgp}{2} \left[ erfc(\frac{h}{\sqrt{2\sgp}})  \exp(\frac{h^2}{2\sgp})\right]^2 \ .
 \NEQ 

A result that goes back to
Levy \cite{Levy, ChoiRoh}, if not earlier,  is 
\begin{theorem}
The joint distribution of the close, $c$, the high, $h$, and the low, $\ell$ is
\BEQ \label{ChoiRoh0}
P(h,\ell,c)=P(B(1)=c, \ell \le \{B(s),\ s \in [0,1]\}\le h) = \sum_{k=-\infty}^{\infty} \phisg(c-2k(h-\ell) ) - \phisg(c-2h- 2k(h-\ell)) 
\NEQ
\BEQ \label{ChoiRoh1}
=  \phisg(c) - \sum_{k=0}^{\infty} \left[\phisg(c-2h- 2k\Delta) + \phisg(c-2\ell+ 2k\Delta) \right]+ \sum_{k=1}^{\infty}\left[\phisg(c-2k\Delta)+ \phisg(c+2k\Delta)\right]
\NEQ
where $\Delta \equiv (h-\ell)$.
\end{theorem}
The symmetric form, \eqref{ChoiRoh1}, not only treats $h$ and $\ell$ symmetrically, but also shows the series is in an alternating form. Here $P(h,\ell,c)$ is a distribution in $h,\ell$ and a density in $c$. 
To calculate the density, we use $p(h,\ell,c) = - \partial_{\ell} \partial_h P(B(t=1)=c, \ell \leq B(t) \leq h)$.
\begin{corollary}
 The density, $p(B(1)=c, \max\{B(s),\ s \in [0,1]\}=h, \min\{B(s)\}=\ell)$, is given by
\BEQ \label{dChoiRoh}
p(h,\ell,c) = \frac{4}{\sigma^2}\sum_{k=-\infty}^{\infty} k^2 a_k(c,\Delta)\phisg(c-2k\Delta) - k(k+1)a_k(c-2h,\Delta) \phisg(c-2h-2k\Delta) \ ,
\NEQ
where $a_k(c,\Delta)\equiv (c-2k \Delta)^2/\sgp-1$.
\end{corollary}
A number of estimators of$\sgg$ given $(h,\ell,c)$ have been proposed \cite{GarmenKlass,McLeish,Meillijson,RogersSatchell}.
The maximum likelihood estimator, $\sggh \equiv argmin_{\sig2} p(h,\ell,c,\sgg) $ was proposed in \cite{BallTorous,Siegmund}.
\begin{corollary}
The density, p($\max\{B(s),\ s \in [0,1]\}=h, \min\{B(s),\ s \in [0,1]\}=\ell)$, satisfies 
  \BEQ \label{densHL}
p(h,\ell) = -  \int_{\ell}^h \partial_{\ell} \partial_hP(B(t=1)=c, \ell \leq B(t) \leq h)dc =
\NEQ \BEQ
\frac{-4}{\sigma^2} \sum_{k=-\infty}^{\infty} k^2 [h_k \phisg(h_k)-\ell_k \phisg(\ell_k)] - k(k+1)[(h_k-2h)\phisg(h_k-2h) - (\ell_k-2h)\phisg(\ell_k-2h)] \ ,
\NEQ
where $h_{k}\equiv h-2k\Delta$ and $\ell_{k}\equiv \ell-2k\Delta$.
\end{corollary}

\section{Density Given High and Close} \label{SheppSect} \label{HCcalcSect}
The classical results in Section 2 are for $t=1$. Our focus for the remainder of the article is on the density
and moments for $t<1$.
We derive the density, $p(B(t)=x |B(1)=c,\ \max\{B(s)\}=h)$ and then compute the first and second moments.
We interpret the high and close as ``statistics'' in the sense of estimation theory.
\begin{theorem} \label{distHCThm}
  The distribution, $F(x,t, h, c) \equiv P(B(t)=x, B(1)=c, B(s) \leq h\ {\rm for}\ s\ \in\ [0,1])$, satisfies 
  \BEQ \label{CHGen}
F(x,t, h, c) \equiv P(B(t)=x, B(s) \leq h\ {\rm for}\ s\ \in\ [0,t]) \times\  P(B'(1-t)=c-x, B'(s) \leq h-x\ {\rm for}\ s\ \in [0,1-t]) \ ,
\NEQ
where $B'$ is a second independent Brownian motion,
\BEQ \label{Pl_CH}
P_{t,x}(h)\equiv P(\max\{B_x(s),\ s \in [0,t]\} \le h) = \phitsg(x) - \phitsg(2h-x) \ ,
\NEQ
\BEQ \label{Pr_CH}
P_{1-t,c-x}(h-x)\equiv P(\max\{B_{c-x}(s),\ s \in [t,1]\} \le h-x) = \phi_{(1-t)\sgp}(c-x) - \phi_{(1-t)\sgp}(2h-c-x) \ .
\NEQ
\end{theorem}

Similar results to  Theorem \ref{distHCThm} for the case of Brownian meanders and excursions appear in \cite{Chung2,DI}, but we have not found precisely this result. One can interpret Theorem \ref{distHCThm}
as a special case of the results in \cite{FPY} where the state space is defined by $B(s) \le h$.
Here $F(x,t,h,c)$ is a distribution in $h$ and a density in $x,c$. 

\begin{corollary} \label{denHCThm}
The conditional density, $p(B(t)=x |h,c) \equiv p(B(t)=x |B(1)=c,\ \max\{B(s)\}=h)$, satisfies:
  \BEQ \label{CHBayes}
p(x,t|h,c)\equiv p(B(t)=x|h,c) = p(B(t)=x, B(1)=c,\ \max\{B(s)\}=h) / p(h,c) \ ,
\NEQ
where the divisor, $p(h,c)$, is given by \eqref{eq:denHC} and 
\BEQ \label{CHDist_Pf0}
p(B(t)=x, B(1)=c,\ \max\{B(s)\}=h) =P_{t,x}(h) q_{1-t,c-x}(h-x)+ q_{t,x}(h) P_{1-t,c-x}(h-x) \ .
\NEQ
Here  
\BEQ \label{Ql_CH}
q_{t,x}(h)\equiv p(\max\{B_x(s),\ s \in [0,t]\} = h) = \frac{d P_{t,x}(h)}{dh} = \frac{2(2h-x)}{t\sgp} \phitsg(2h-x)
\NEQ
\BEQ \label{Qr_CH}
q_{1-t,c-x}(h-x) \equiv p(\max\{B_{c-x}(s),\ s \in [t,1]\} = h-x) = \frac{2(2h-c-x)}{(1-t)\sgp} \phi_{(1-t)\sgp}(2h-c-x) \ .
\NEQ
\end{corollary}

Equation \eqref{CHDist_Pf0} simply states that if the realization goes through the points $(t,x)$ and $(1,c)$ and has a high value of $h$,
then either it reaches $h$ in $[0,t]$ or in $[t,1]$. Equation \eqref{CHDist_Pf0} is the kernel of the Chapman-Kolmogorov representations for this restricted Brownian motion problem.
\begin{lemma}
$F(x,t, h, c)$ is the difference of four Gaussians:
\BEQ
\label{CHMomGenSum}
F(x,t, h, c)  = f_1(x,t,c) -f_2(x,t,c,h) - f_3(x,t,c,h) +f_4(x,t,c,h) \ . 
\NEQ
The  $f_i$ are of the form: 
\BEQ \label{fiDef}
f_i(x,t,h,c)=\frac{1}{2\pi\sig\sig_t}\exp\left(-\frac{(1-t)(x-a_i)^2 +(x-b_i)^2 t}{2t(1-t)\sgp}\right)=
\phi_{\sgp_t}(x- \mu_i(t,c,h)) )\frac{\exp^{-g_i(c,h)}}{\sqrt{2\pi}\sig}
\NEQ
where $\sigma_t^2 \equiv t(1-t)\sigma^2$, 
$a_1=0$, $b_1=c$, $a_2=0$, $b_2=2h-c$, $a_3=2h$, $b_3=c$, $a_4=2h$, $b_4=2h-c$. In \eqref{fiDef}, $\mu_i(t,c,h)$ and $g_i$ are defined by
\BEQ \label{muDef}
\mu_i(t,c,h) \equiv a_i(1-t) +b_it\ ;\quad g_i \equiv\frac{[a_i^2(1-t) +t b_i^2]-\mu_i^2}{2\sigma_t^2}=\frac{(a_i - b_i)^2}{2\sigma^2} \ .
\NEQ
Thus,
$\mu_1= ct$, $g_1=c^2/2\sgp$, $\mu_2= (2h-c)t$, $g_2=(2h-c)^2/2\sgp$, $\mu_3=2h(1-t) +ct$, $g_3=(2h-c)^2/2\sgp$ and $\mu_4= 2h-ct$, $g_4=c^2/2\sgp$.
\end{lemma}

Note that $f_i(x=h,t,h,c)=  \phi_{t\sgp}(h)\phi_{(1-t)\sgp}(h-c)$. 
The equality of the four terms at $x=h$ will allow us to cancel terms when we integrate by parts.
We also define $\psi_i = \exp^{-g_i(c,h)} / \sqrt{2\pi}\sig$ so $\psi_1=\psi_4=\phi_{\sgg}(c)$ and $\psi_2=\psi_3=\phi_{\sgg}(2h-c)$.
To simplify our calculations, observe
\BEQ \label{CHder}
\partial_h f_i = \left[\frac{(x-\mu_i)}{\sigma_t^2} \partial_h \mu_i -\partial_h g_i\right] f_i \quad {\rm and}\ \partial_x f_i = -\frac{(x-\mu_i)}{\sigma_t^2} f_i \ .
\NEQ
Note $f_1$ is independent of $h$ and therefore may be ignored. Derivatives of 
$F(x,t, h, c)$ with respect to $h$ only enter through $h$ dependencies in $\mu_i$ and $g_i$.
We now evaluate the moments with respect to $x$ for a given time, $t$, and fixed $(h,c)$.

\begin{theorem} \label{momThm}
Consider $M_m(t, h, c) \equiv \int_{\infty}^h x^m p(x,t, h, c)dx$. The zeroth moment is  $M_0(t,c,h) =p(h,c)$ where $p(h,c)$ is given in \eqref{eq:denHC}.
\BEQ
\label{CHMom1}
M_1\equiv \phi_{\sgg}(c)[1 +erf(\frac{ct-h}{\sqrt{2}\sigma_t})] +\phi_{\sgg}(r)[\frac{2hr}{\sgg}-1 +p_{h,r,t} erf(\frac{h-rt}{\sqrt{2}\sigma_t})] - 4 r t(1-t) \phi_{\sgg}(r) \phi_{\sigma_t^2}(h-rt) 
\NEQ
where $r\equiv 2h-c$ and $p_{h,r,t}\equiv (1-2t) +\frac{2r(rt-h)}{\sgg}$.
\BEQ
\label{CHMom2}
M_2=2(2h-ct) \phi_{\sgg}(c)[1 +erf(\frac{ct-h}{\sqrt{2}\sigma_t})] +
\NEQ
\BEQ
2\phi_{\sgg}(r) \left(\left[r t(1-t) + q_1(h,c,t) + q_2(h,c,t) erf(\frac{h-rt}{\sqrt{2} \sigma_t})\right]
  - 4 r h t(1-t) \phi_{\sigma_t}(h-rt) \right) \ ,
\NEQ
where  $q_1(h,c,t)\equiv -(rt^2 +(1-t)(2h-rt) ) +r (h^2+ (h-rt)^2)/\sgg$ and $q_2(h,c,t)\equiv (2h(1-t)-rt) + 2hr(rt-h)/\sgg$
\end{theorem}
Note $M_1(t=1,c,h) =c p(c,h)$ and $M_2(t=1,c,h) =c^2 p(c,h)$ as must be.


To compute the moments, we use 
\BEQ
\label{CHMomi}
M_m(h,c)=\int_{\infty}^h x^m \partial_h F(x,t, h, c)dx  = \sum_{i=2}^4 s_i \int_{-\infty}^h x^m \left[-\partial_h \mu_i \partial_x f_i -\partial_h g_i f_i\right]  =
\NEQ
\BEQ
\label{CHMomiB}
\sum_{i=2}^4 s_i \int_{-\infty}^{h}\left[ m x^{m-1} \partial_h \mu_i -\partial_h g_i x^m \right]f_idx =
\NEQ
\BEQ
\label{CHMomiC}
 \sum_{i=2}^4 s_i \psi_i\int_{-\infty}^{h-\mu_i}
\left[m(x+\mu_i)^{m-1} \tau_i - \frac{2(2h-c)}{\sigma^2}(x+\mu_i)^m(1-\delta_{i,4})\right] \phi_{\sigma_t^2}(x)dx \ .
\NEQ
Here  $s_i = -1$ for $i=2,3$ and $s_i=1$ for $i=1,4$ and
we define $\tau_i \equiv \partial_h \mu_i$ so that $\tau_2= 2t$, $\tau_3= 2(1-t)$, $\tau_4= 2$. 
To go from \eqref{CHMomi} to \eqref{CHMomiB}, we use that the three terms evaluated at $x=h$ cancel.
The remainder of the evaluation of the moments $M_1(t,h,c)$ and $M_2(t,h,c)$  is deferred to the Appendix.
\qed

Given the moments, $M_i(t,h,c)$, the conditional mean and conditional second moment are  $E[B(t|h,c)]=M_1(t,h,c)/p(h,c)$
and $E[B^2(t|h,c)]=M_2(t,h,c)/p(h,c)$. We treat $E[B(t|h,c)]=M_1(t,h,c)/p(h,c)$ as an estimator of $B(t|h,c)$. An alternative
estimator is the maximum likelihood estimate given by maximizing \eqref{CHDist_Pf0} with respect to $x$ for each time $t$.



\section{Density Given High, Low and Close} \label{HLCcalcSect}

We now consider using the open, high, low and close together as statistics to estimate a realization of a Brownian process.
After writing down the density conditional on these statistics, we spend the bulk of this section evaluating
the moments of the density as summarized by \eqref{momThmHLC}.
We begin by applying Chapman-Kolmogorov equation to $P(B(t)=x, \ell \leq B(s) \leq h)$:
\begin{theorem} \label{distHLCThm}
  Let $G(x,t, h, \ell, c) \equiv P(B(t)=x, B(1)=c, \ell \leq B(s) \leq h | {\rm for}\ s \in [0,1] )$,
  $Q(x,t, h,\ell) \equiv P(B(t)=x, \ell \leq B(s) \leq h|s\le t)$ and 
  $Q_R(x,t, h,\ell,c) =  P(B'(1-t)=c-x, \ell-x \leq B'(s) \leq h-x\ |s\le 1-t)$. Here $B'$ is a second independent
Brownian motion. Then
\BEQ \label{ChoiRoh}
P(B(t)=x, \ell \leq B(s) \leq h|s\le t) = \sum_{j=-\infty}^{\infty} \bigg{[} \phitsg(x-2j(h-\ell)) - \phitsg(x-2h+ 2j(h-\ell)) \bigg{]} \ , 
\NEQ
\BEQ \label{ChoiRohR}
P(B'(1-t)=c-x, l-x \leq B'(s) \leq h-x) = \sum_{k=-\infty}^{\infty} \bigg{[} \phihsg(x-c+2k\Delta) - \phihsg(x-(2h-c) -2k\Delta) \bigg{]} 
\NEQ
where $\Delta$ is the high - low on $[0,1]$, $\Delta\equiv h-\ell$ and $\hat{t}\equiv (1-t)$.
The probability distribution, $P(B(t)=x, \ell \leq B(s) \leq h)$, satisfies
\BEQ \label{CHLDef}
G(x,t, h, \ell, c) = Q(x,t, h,\ell) Q_R(x,t, h,\ell,c)
\NEQ
\end{theorem}

{\it Proof:} We apply \eqref{ChoiRoh0} in the time interval $s \in [0,t]$ to yield \eqref{ChoiRoh} and to $s \in [t,1]$ to yield \eqref{ChoiRohR}.
The Markovian property yields \eqref{CHLDef}. \qed

One can interpret Theorem \ref{distHLCThm} as a special case of the results in \cite{FPY} where the state space
is defined by $\ell \le B(s) \le h$.
Clearly, $Q_R(x,t, h,\ell, c)= Q(c-x,1-t, h-x,\ell-x)$.
Here $G(x,t,h,\ell,c)$ is a distribution in $h,\ell$ and a density in $x,c$. 
To derive the density of $p(x;t,h,\ell,c)$, we need to consider four terms, the probability that both the high and low
are to the left of $t$, the probability that just the low is to the left of $t$, the probability that just the high is to
the right of $t$ and the probability that both the high and the low are to the right of $t$.
\begin{corollary}
The density, $p(x;t,h,\ell,c)\equiv  p(B(t)=x |B(1)=c,\ \max\{B(s)\}=h,\ \min\{B(s)\}=\ell)$, satisfies
\BEQ
\label{CHLGen2}
p(x;t,h,\ell,c) = - \partial_{\ell}\partial_h G(x,t, h,\ell, c) = -\partial_{\ell}\partial_h  Q(x,t, h,\ell) Q_R(x,t, h,\ell,c) \ .
\NEQ
Furthermore,
\BEQ 
p(x;t,h,\ell,c) =p(t_{\ell}<t, t_{h}<t)+p(t_{\ell}<t, t_{h} \ge t)+ p(t_{\ell} \ge t, t_{h} < t) +  p(t_{\ell}\ge t, t_{h} \ge t) 
  \label{CHLDefb}
\NEQ
where $t_{\ell}$ is the first time that $B$ reaches its minimum and $t_{h}$ is the first time that $B$ reaches its maximum.
\end{corollary}
Analogous to \eqref{CHDist_Pf0}, equation \eqref{CHLGen2} is the kernel of the Chapman-Kolmogorov representation for this restricted Brownian motion problem. The four terms in \eqref{CHLDefb} correspond to applying the product rule of calculus to $Q(x,t, h,\ell) Q_R(x,t, h,\ell,c)$.
As in \eqref{fiDef}, the generator, $G(x,t, h, \ell, c)$, is composed of a sum of Gaussians in $x$. The remainder of this section and the Appendix
are devoted to evaluating the moments of $p(x;t,h,\ell,c)$.
\begin{theorem} \label{Gexpand}
  The probability distribution $G(x,t, h, \ell, c)$ has the following representations:
\BEQ \label{CHLSum1}
G(x,t, h,\ell, c) \equiv \sum_{j,k>-\infty}^{\infty} \sum_{i=1}^4  \frac{s_{i}}{2\pi\sig \sig_t} \exp\left(\frac{-(x-a_{ij})^2}{2t\sgp} -\frac{(x-b_{ik})^2 }{2(1-t)\sgp}\right)
=\ \sum_{j,k>-\infty}^{\infty}\sum_{i=1}^4  s_{i} f_{ijk}(x) \ . 
\NEQ
Here  $s_i = -1$ for $i=2,3$ and $s_i=1$ for $i=1,4$ and $\sigma_t^2\equiv t(1-t)\sigma^2$.
The parameters are defined as $a_{1,j}=2j \Delta$, $b_{1,k}=c-2k\Delta$, $a_{2,j}=2j \Delta$,
$b_{2,k}=2h-c+ 2k\Delta$,$a_{3,j}=2h-2j\Delta$, $b_{3,k}=c-2k\Delta$, $a_{4,j}=2h-2j\Delta$, $b_{4,k}=2h-c+2k\Delta$.
The shifted Gaussian representation is
\BEQ \label{CHLSum2}
G(x,t, h,\ell, c)=\sum_{j,k>-\infty}^{\infty} \sum_{i=1}^4 \frac{s_{i}}{\sqrt{2\pi}\sigma}\phi_{\sigma_t^2}(x- \mu_{ijk}(h,\ell)) e^{-g_{ijk}(h,\ell)}
= \sum_{j,k>-\infty}^{\infty} \sum_{i=1}^4 s_{i}\phi_{\sigma_t^2}(x-\mu_{ijk}(h,\ell)) \psi_{ijk}(h,\ell)
, \
\NEQ
where $\mu_{ijk}(t,c,h) \equiv (a_{ij}(1-t) +b_{ik}t)$, 
and $g_{ijk} =\left([a_{i,j}^2(1-t) +t b_{i,k}^2]-\mu_{ijk}^2\right)$ $/(2\sigma_t^2)=(a_{i,j} - b_{i,k})^2/2 \sigma^2 $
and $\psi_{ijk}(h,\ell) \equiv  e^{-g_{ijk}(h,\ell)} /\sqrt{2\pi}\sigma$.

Let $v_{j,k} =2(j(1-t) +kt)$, $\vt_{j,k} =2(j(1-t) -kt)$, $w_{j,k} =2(j+k)$, $\wt_{j,k} =2(j-k)$.
Then
$\mu_{1,j,k}= ct+\vt_{j,k}\Delta$, $g_1= (c-w_{j,k}\Delta)^2/2\sgg$,
$\mu_2= (2h-c)t+v_{j,k}\Delta$, $g_2=(2h-c-\wt_{j,k}\Delta)^2/2\sgg$, $\mu_3=2h(1-t) +ct-v_{j,k}\Delta$, $g_3=(2h-c-\wt_{j,k}\Delta)^2/2\sgg$ and $\mu_4= 2h-ct-\vt_{j,k}\Delta$, $g_4=(c-w_{j,k}\Delta)^2/2\sgg$.
\end{theorem}

As in Section \ref{HCcalcSect}, we evaluate the moments in $x$ for a given time, $t$, and fixed $(h,\ell,c)$.
\begin{lemma}
  The moments, $M_m(t,h,\ell,c)$,
\BEQ  \label{CHLMomGen}
M_m(t,h,\ell,c) \equiv -\int_{\ell}^h x^m \partial_h \partial_{\ell} G(x,t, h,\ell, c)=\int_{\ell}^h x^m \sum_{j,k>-\infty}^{\infty} \sum_{i=1}^4 -s_{i} \partial_h \partial_{\ell} f_{ijk}(x,h,\ell) \ \ .
\NEQ
The $(i,j,k)$th term inside the integral satisfies
\BEQ \label{CHLMomijk}
-\partial_h \partial_{\ell} \phi_{\sigma^2_t}(x- \mu_{ijk}(t,h,\ell))e^{-g_{ijk}(h,\ell,t)} = H_{ijk}(x- \mu_{ijk}) \phi_{\sigma_t^2}(x- \mu_{ijk}(t,h,\ell))e^{-g_{ijk}} \ \ .
\NEQ
Here $H_{ijk}(z)$ is a quadratic polynomial in $z$ defined as
\BEQ \label{CHLMomijk2}
H_{ijk}(z)\equiv -[\frac{z \tau_{ijk} }{\sigma_t^2} -\partial_h g_{ijk}] [\frac{z \tauh_{ijk}}{\sigma_t^2}  - \partial_{\ell} g_{ijk}] + \frac{\tau_{ijk}\tauh_{ijk}}{\sigma_t^{2}} + \partial_{\ell}\partial_h g_{ijk}\ .
\NEQ
where $\tau_{ijk} \equiv \partial_{h} \mu_{ijk}$ and $\tauh_{ijk} \equiv \partial_{\ell} \mu_{ijk}$. Thus $\tau_{1jk} =\vt_{j,k}=-\tauh_{1jk}$, $\tau_{2jk} =2t+v_{j,k}$, $\tauh_{2jk} =-v_{j,k}$, $\tau_{3jk}=2(1-t)-v_{j,k}$, $\tauh_{3jk} =v_{j,k}$, $\tau_{4jk} =2-\vt_{j,k}$ and $\tauh_{4jk} =\vt_{j,k}$.  
Here $\tau_{.}$ and $\tauh_{.}$ have no dependence on $h$ and $\ell$. 
\end{lemma}

We group the terms in \eqref{CHLMomijk2} by powers of $z$ and define
$A_{ijk} =  \tau_{ijk}  \tauh_{ijk}$, $B_{ijk} = [\tau_{ijk}\partial_{\ell} g_{ijk}+\tauh_{ijk}\partial_{h} g_{ijk}]$ and
$C_{ijk} =\Gamma_{ijk}+\sigma_t^{-2} \tau_{ijk}  \tauh_{ijk} $ and $\Gamma_{ijk}\equiv -\partial_h g_{ijk}\partial_{\ell} g_{ijk} + \partial_{\ell}\partial_h g_{ijk}$.
Note that $\Gamma_{4jk} = \Gamma_{1jk} =  (2*g_{1jk}-1)w_{jk}^2/\sgg$ and $\Gamma_{3jk} = \Gamma_{2jk} =  (2*g_{2jk}-1)\wt_{jk}(\wt -2)/\sgg$.
Thus $H_{ijk}(z) = -A_{ijk}z^2/\sigma_t^4 +B_{ijk}z/\sigma_t^2 +C_{ijk}$.

To simplify the moment calculation, we evaluate the derivatives by $h$ and $\ell$ and recast them as derivatives with respect to $x$ so that we can integrate by parts:
\BEQ \label{CHLDhlvsDxl}
\partial_h \partial_{\ell} f_{ijk} = \tau_{ijk}\tauh_{ijk}\partial_x^2  f_{ijk} + B_{ijk}\partial_x  f_{ijk} - \Gamma_{ijk} f_{ijk} \ \ .
\NEQ
Note that $\sum_{i=1}^4 s_i f_{ijk}(x=h)=0$,
$\sum_{i=1}^4 s_i \partial_h f_{ijk}(x=h)=0$ and $\sum_{i=1}^4 s_i \partial_l f_{ijk}(x=h)=0$. This allows us to integrate by parts and drop terms.

We define the moments, $M_m$, where the limits of integration, $H$ and $L$, are to be set to $h$ and $\ell$ after
we differentiate $\partial_h \partial_{\ell} G$. This is done because integration by parts should not include the
dependence on the limits of integrations.

We integrate by parts and find from \eqref{CHLDhlvsDxl}:
\BEQ \label{CHLMomGenM1} 
M_m(t,h,\ell,c)=  \sum_{j,k>-\infty}^{\infty} \sum_{i=1}^4 s_{i} \int_L^H mx^{m-1} \left[A_{ijk}\partial_{x}f_{ijk} 
  +B_{ijk} f_{ijk} \right]   + x^m \Gamma_{ijk}  f_{ijk}
\NEQ
\Beq
+ \sum_{i,j,k} s_i (\Bc_{ijk}(x=h) -\Bc_{ijk}(x=\ell) )
\Neq
where $\Bc_{ijk}$ is the boundary term. In the Appendix \ref{AppHLCBC}, we show that the boundary terms sum to zero. 

 In this section, we will often need the triple sum, $\sum_{j,k>-\infty}^{\infty} \sum_{i=1}^4$. For notational simplicity, we replace the triple sum with a simple $\sum_{ijk}$ when appropriate.

\begin{lemma}
The moments have the representation:
\BEQ \label{HLCmom0}
M_m(t, h,\ell,c) = \sum_{ijk} s_{i} \psi_{ijk}\left[\frac{-m A_{ijk}}{\sgg_t}G_{m-1,1}(\mu_{ijk}) +m B_{ijk}G_{m-1,0}(\mu_{ijk},\sgg_t)+\Gamma_{ijk}G_{m0}(\mu_{ijk},\sgg_t)   \right] \ .
\NEQ
where the 
\BEQ
\label{Gmn}
G_{mn}(\mu,h,\ell,\sigma)\equiv\int_{\ell-\mu}^{h-\mu} (x+\mu)^m x^n \phi_{\sigma^2}(x) \ .
\NEQ
As in Theorem \ref{Gexpand}, $\psi_{ijk} \equiv e^{-g_{ijk}(h,\ell,t)}  /\sqrt{2\pi} \sigma$.
\end{lemma}

\Prf We substitute the definitions in \eqref{Gmn} into \eqref{CHLMomGenM1}.


In Appendix \ref{GmnEval}, we evaluate the functions $G_{mn}()$ in terms of $\phi_{\sggt}(h-\mu_{ijk})$, $\phi_{\sggt}(\ell-\mu_{ijk})$ and the
corresponding error functions. Collecting terms from above and using the Appendix \ref{GmnEval} yields

\begin{theorem} \label{momThmHLC}
For $m \le 2$, equation \eqref{HLCmom0} becomes
\BEQ \label{HLCM1g}
M_m(t)=\sum_{ijk} s_{i}\psi_{ijk}\left[a_{ijk}^{(m)} \phi_{\sggt}(h-\mu_{ijk}) -\ah_{ijk}^{(m)} \phi_{\sggt}(\ell-\mu_{ijk})
+  e_{ijk}^{(m)}R_{\sigma_t}(h,\ell,\mu_{ijk})\right] \ ,
\NEQ
where $R_{\sigma_t}(h,\ell,\mu_{ijk})\equiv [ E_{\sigma_t}(h-\mu_{ijk}) - E_{\sigma_t}(\ell-\mu_{ijk})]$ and $E_{\sigma}(x)$ is
the scaled $erf$ function, $E_{\sigma}(x)\equiv .5* erf(x/\sqrt{2}\sigma)$.
For $m=1$, the coefficients are
\BEQ \label{HLCMcf1}
a_{ijk}^{(1)}= \ah_{ijk}^{(1)}= A_{ijk}- \Gamma_{ijk}*\sgg_t \ \ \ , \ \ \
e_{ijk}^{(1)}=B_{ijk}+ \Gamma_{ijk}\mu_{ijk}  \ .
\NEQ
For $m=2$, the coefficients are
\begin{align}
a_{ijk}^{(2)}  &=2h A_{ijk} - 2B_{ijk}\sgg_t - \Gam_{ijk}\sgg_t (\mu_{ijk}+h)  \ \ , \nonumber \\  
\ah_{ijk}^{(2)}&=2\ell A_{ijk} -2B_{ijk}\sgg_t - \Gam_{ijk}\sgg_t(\mu_{ijk}+\ell) \ \ , \nonumber \\ 
e_{ijk}^{(2)}  &=2B_{ijk}\mu_{ijk}-2A_{ijk} + \Gam_{ijk}*(\mu_{ijk}^2+\sgg_t) \ \ .
\label{HLCMcf2}
\end{align}
For $m=0$, $a_{ijk}^{(2)}=0$, $\ah_{ijk}^{(2)}=0$ and $e_{ijk}^{(2)}=\Gam_{ijk}$.
\end{theorem}

\Prf We substitute in the $G_{mn}$ expressions into \eqref{HLCmom0} and collect terms.

Some further simplifications of the coefficients in \eqref{HLCMcf1}-\eqref{HLCMcf2} can be found in Appendix \ref{AppSimp} for $m\le 2$. 
When $m>2$, the terms multiplying $E_{\sigma}(h-\mu_{ijk})$ and $E_{\sigma}(\ell-\mu_{ijk})$ are different.
\begin{corollary} \label{corM0}
$M_0(t,h,\ell,c)= p(h,\ell,c)$ as given by \eqref{dChoiRoh}.
\end{corollary}

\Prf See the Appendix \ref{AppM0}. 

We treat $E[B(t|h,\ell,c)]=M_1(t,h,\ell,c)/p(h,c)$ as an estimator of $B(t|h,\ell,c)$. 
To evaluate \eqref{HLCmom0} numerically, we need to truncate the expansion in $j$ and $k$. Luckily, the Feller distribution of \ref{FellerSect} shows that very few realizations have small values of $\Delta$.
Thus the double expansion for $j$ and $k$ converges quickly for the vast majority of the Brownian realizations.

A second method to evaluate the probability of \eqref{CHLDef} is 
to numerically evaluate $Q$, $\partial_h Q$, $\partial_{\ell} Q$
and $\partial_{\ell}\partial_h Q$ in \eqref{ChoiRoh}
and to numerically evaluate $Q_R$, $\partial_h Q_R$, $\partial_{\ell} Q_R$ and $\partial_{\ell}\partial_h Q_R$ in \eqref{ChoiRohR}.
We then numerically integrate the moments of \eqref{probQQr}.
\BEQ \label{probQQr}
prob(x,t,h,\ell,c)=Q(x,t,h,\ell) \partial_{\ell}\partial_h Q_R + \partial_h Q(x,t,h,\ell) \partial_{\ell} Q_R + \partial_{\ell} Q \partial_h Q_R+ Q \partial_{\ell}\partial_h Q_R(x,t,h,\ell,c)
\NEQ
times $x^m$ from $x=\ell$ to $x=h$.
Each of the terms in the integral involves truncating only in one of $j$ or $k$. Thus the additional work involved in evaluating $Q$ and $Q_R$ at many points to evaluate the integral is partially compensated by the single infinite sums as opposed to a doubly infinite sum. 

An alternative estimator is the maximum likelihood estimate given by maximizing likelihood of $p(x,t|h,\ell,c)$ with respect to $x$ for each time $t$.
Here $p(x,t|h,\ell,c) = -\partial_h \partial_{\ell}G(x,t, h, \ell, c) / p(h,\ell,c)$. Using the series representation yields 
\BEQ \label{pcondHLC}
p(x,t|h,\ell,c) =\frac{\sum_{j,k>-\infty}^{\infty} \sum_{i=1}^4 \frac{s_{i}}{\sqrt{2\pi}\sigma} H_{ijk}(x- \mu_{ijk}(h,\ell)) \phi_{\sigma_t^2}(x- \mu_{ijk}(h,\ell)) e^{-g_{ijk}(h,\ell)}}{p(h,\ell,c)}
\NEQ
In practice, the estimator $E[B(t|h,\ell,c)]=M_1(t,h,\ell,c)/p(h,c)$ is much faster to evaluate than the maximum likelihood estimate.

\section{Numerical Methods} \label{NumSect}

Simply put, we generate a large number of Brownian paths, bin the paths in $(close,\ max,\ min)$ space
and calculate the mean and variance for each time and bin. We order the coordinates of phase
space, $(q_1,q_2, q_3)$, so that $q_1 = B(1)$, $q_2=\ max_{0 \le t \le 1}B(t)$ and $q_1=\ min_{0 \le t \le 1}B(t)$.
We also consider the case where we replace one or more of these operators with $argmax$ or $argmin$.
The results for the $argmax$ case are found in \cite{RiedelArgMx}.

A very straightforward algorithm is

1) Specify a timestep, $dt$, a number of bins in each direction $nbins$, and a number of sample paths, $N_{samp}$
with typically $N_{samp} \approx \kappa\ {\rm nbins}^3$ where $\kappa$ denotes the typical number of simulations in a bin.
More generally, for any choice of grids for the bins, we want at least $\kappa$ simulations in each bin where $\kappa$ is a large number.
Generate a large array of scaled Gaussian random variables, size $(N_{samp}, 1/dt)$.
Cumsum them to generate an array of Brownian paths. We often use a nonuniform time step where the time step is smaller near $t=0$ and near $t=1$.

2) In the first phase space direction, compute bin boundaries so that the number of curves are roughly equal in each bin.
For each one dimensional bin, compute bin boundaries in the second coordinate direction so that the number of bins
is roughly equal. Finally, for each of the two dimensional bins, compute bins in the third direction.

3) For each of the $nbin^3$ bins, assign a triple index, $\left(i,j,k\right)$ bins,
compute the mean of the coordinates, $(\qb_1,\qb_2,\qb_3)$, and compute the mean, $\mu(t;\qb_1,\qb_2,\qb_3)$,
and variance, $V(t;\qb_1,\qb_2,\qb_3)$, of $\{B(t)\}$ in the bin.

4) Test for convergence of $\mu(t;\qb_1,\qb_2,\qb_3)$ and $V(t;\qb_1,\qb_2,\qb_3)$ in $N_{samp}$, $nbins$, and $dt$.
This involves interpolation as grids boundaries are random functions of the particular ensemble of paths.
Note that the grid boundaries for the first coordinate direction are independent of the second two coordinate directions
but that the average value of $q_1$ will depend on all three indices, $(i,j,k)$. We find that interpolation from
one grid to another grid broadens the width of the peaked functions especially when $argmax$ is one of the given variables.

There is a bias versus variance tradeoff. If the bins are too large, the variation of the mean and variance will be obscured.
If the bins are too small, there will be too few curves in each bin and the sample variance will dominate.
Each of the close, max and min have a Gaussian or half Gaussian distribution individually so the tails of the distribution will be spread out.
The situation is actually somewhat better as the high and low are exponentially distributed given the closing value.
Nevertheless, exponential distributions have very few points in the tail of the distribution. Again, a low density of
curves will significantly inflate the size of the tail bins and thereby add larger bias to the the computation of
the bin variance. Thus convergence of the mean and variance on the outermost bins is tenuous.
When we compute population average variance, we are tempted to downweight or even exclude the outer bins. While
this is probably a smart thing to do, we report the simple ensemble average instead of a more complex limit reducing the underweighting as the bin size goes to zero.

Assume that we wish to generate bins in the $\qb$ direction. We sort the Brownian realization in the $\qb_1$ direction.
To generate the grid boundaries, we initially tried equi-spaced quantile bins. This results in very large bins in the low
density region. These large bins result in bias to our estimates for both the expectation and variance estimates.
Let the density of points/curves be $n(\qb)$. To reduce the the size of the largest bins, we select bin boundaries to
keep $\int_{\qb_k}^{\qb_{k+1}} n(\qb)^\alpha d\qb$ to be approximately equal where $\{\qb_k\}$ are the bin boundaries. We
use $\alpha=.7-.75$ while $\alpha=1$ corresponds to equal quantiles.
We find that first and last bins converge much very slowly in $(nSim, nbin)$ space especially using a quantile
based gridding. Using equal bins of $n(\qb)^\alpha$ partially but not completely alleviates this problem.


Wiener's Fourier representation of Brownian paths on $[0,1]$ is
\BEQ \label{BFourier}
B(t) = \xi_0 t + \sum_{n=1}^{\infty} \xi_n \frac{sin(nt)}{\pi n},\ \ {\rm where\ } \{\xi_k \} {\rm\ are\ independent\ normal.}
\NEQ
Given an ensemble of Brownian paths, $\{B_i(t)\}$, we can create an equivalent ensemble of Brownian paths, $\{B_i(t,c)\}$, with right endpoint $c$, using the formula:
$B_i(t,c) \equiv B_i(t) - (B_i(t=1) -c) t$.
This allows us to take one set of Brownian paths and use them on a grid of final values. This significantly reduces the number
of realizations we need to cover phase space. Thus if the closing value is the first parameter direction that we examine, a
3-dimensional parameterization is reduced to a sequence of two-dimensional parameterizations.



\section{Single Conditional Value \label{SCVsect} }
\subsection{Brownian Bridge}
\begin{figure}[htbp]
\centering
\includegraphics[scale=.5]{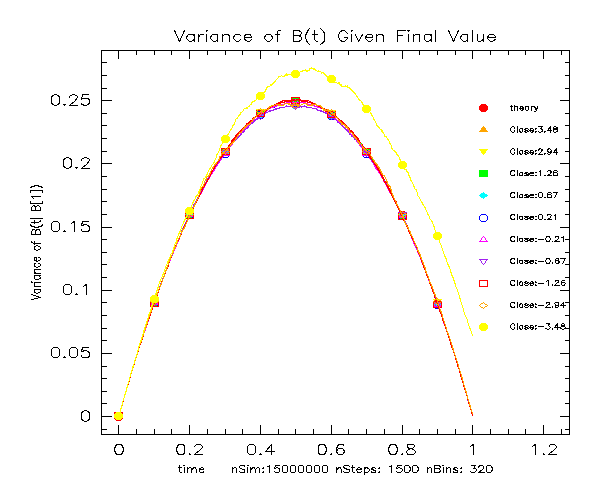}
\caption{$Var[B(t|c)]$ for various final values, $c$ }
\label{fi:VarArrGivenClose}
\end{figure}

We begin with plots of our simulation for the Brownian bridge case, i.e. Brownian motion constrained to a given closing value.
For this simulation, we use 15 million simulations with nsteps=1500. For a given value of $B(1)=c$, the simulation yields a straight line
in time for $E[B(t)|B(1)=c]$.
Figure \ref{fi:VarArrGivenClose} plots the time dependent variance,  $Var[B(t)|B(1)=c]$ for a variety of $c$.
The closing values are chosen to be the values inbins number, $(0,2, \ldots nbin-3, nbin-1)$
where the third through eigth bin are equi-spaced in bin number.
The theoretical value is $t(1-t)$ and is displayed as the red curve in Figure \ref{fi:VarArrGivenClose}.
All but the first and last curve match the theoretical values. 
This occurs because the first and last bins cover a very large range of $c$. We are averaging different values of $E[B(t)|B(1)= close]$
and the squared bias is miscounted as variance. 

\subsection{Given High}

To calculate the probability, $p(x,t, h)$, we integrate $p(x,t, h)= \int_{-\infty}^h p(x,t, h, c)dc$
using\eqref{CHDist_Pf0}-\eqref{Qr_CH}. The result is
\begin{theorem}
The probability density, $p(x, t, h) \equiv p(B(t)=x |\max\{B(s)\}=h)$ satisfies
  \BEQ \label{Prbxth}
p(x,t, h) = 2[\phitsg(x) - \phitsg(2h-x)]\phi_{(1-t)\sigma^2}(h-x)  +  \frac{2(2h-x)}{t\sigma^2} \phitsg(2h-x) erf(\frac{h-x}{\sqrt{2(1-t)\sigma_t}}).
\NEQ
\end{theorem}
The theoretical values of  $E[B(t)|h]$ and  $Var[B(t)|h]$ can be calculated by computing moments with respect to \eqref{Prbxth} and
then dividing by $p(h) = 2 *\phi_{\sigma^2}(h)$, $h\ge0$.
Unfortunately, we have not found a tractable analytic form from the integrals and therefore we compute them numerically \cite{FA}.
Another, very abstract, expresion for $p(x, t| h)$ can be found in \cite{BertoinPitmanChavez}.

\begin{figure}[htbp]
\centering
\includegraphics[scale=.5]{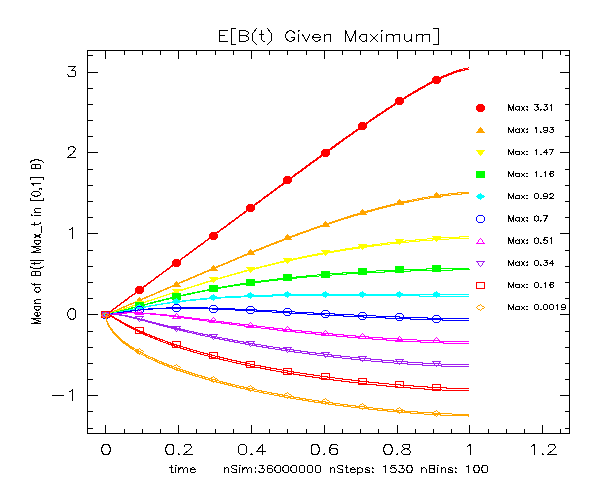}
\caption{Expectation of $B(t)$ given $max\{B(s)\}=h$ for various values of the high, $h$. If the maximum is small, the expectation decreases
for most of its domain. If the maximum is large, the expectation increases. }
\label{fi:MeanArrGivenMax}
\end{figure}

\begin{figure}[htbp]
\centering
\includegraphics[scale=\FSCL]{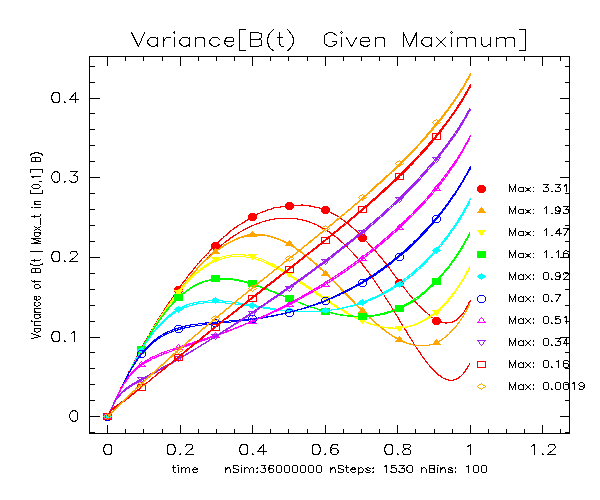}
\caption{Variance of  $B(t)$ given $max\{B(s)\}=h$ for various values of the high, $h$.}
\label{fi:VarArrGivenMax}
\end{figure}

Figure \ref{fi:MeanArrGivenMax} plots the expectation of $B(t)$ for ten values of the high.
Not surprisingly, if the high occurs near $t=0$, the expectation
decreases monotonically beyond the argmax, $\theta$, and decreases faster for smaller $t$. Let $f(t,h)\equiv E[B(t)|\max B=h]$. It appears
that $f$ is smooth in $t$ and $|\frac{\partial f}{\partial t}|$ is decreasing in time. For large values of $h$, $f(t,h)$
grows approximately linearly. We see that the zero of $f(t=1,h)$ occurs somewhere between $.68<h<.95$. Using \eqref{eq:EVCH}, we see the
precise value is .7517915247. Figure \ref{fi:VarArrGivenMax} plots the variance of a bin as a function of time.
Again, the computed variance includes the squared bias from effectively assuming that expectation is constant in each bin.
Since $f(t,h)$ varies from the smallest
value of $h$ in the bin to the largest value of $h$ in the bin, we are systematically overestimating the variance.
For this particular computation, we define $vrAvg$ to be the time and ensemble average of the variance.
Looking at the dependence as a function of $nbins$, the number of bins, we find $vrAvg(nbins=80)=0.16033$,
$vrAvg(nbins=160)=0.16021$, $vrAvg(nbins=320)=0.16018$ and $vrAvg(nbins=480)=0.16017$.
Knowing the value of the high is slightly better at reducing the time averaged variance since $vrAvg< 1/6$.

Returning to Figure \ref{fi:VarArrGivenMax}, we see that that $var(t,h)\equiv Var[B(t)|\max B= h]$ is monotonically increasing for small values of $h$, up to at least $h=.67$. For larger values of $h$, the variance is non-monotone. This non-monotonicity occurs because at large values of $h$ , the maximum of $B(t)$ is likely to be near $t=1$.
In these simulations, we use an ensemble of 36,000,000 realizations computed with 1530 steps and bin the results into 100 bins.

\begin{figure}[htbp]
\centering
\includegraphics[scale=.5]{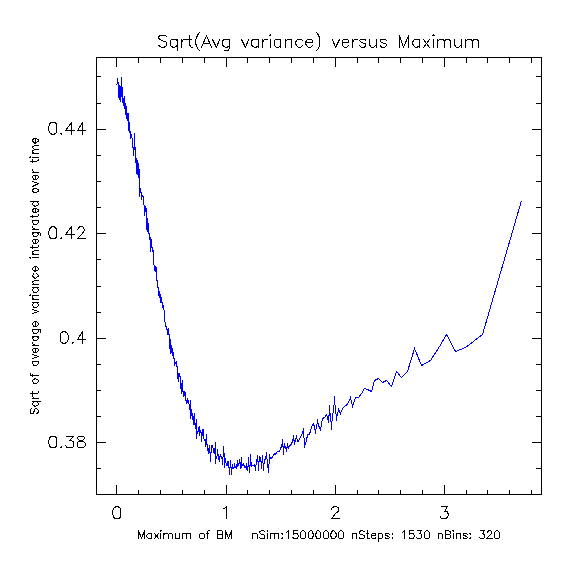}
\caption{Square root of time average of  $Var[B(t)|h]$ versus the high, $h$. We strongly suspect that the final point on the curve
is incorrect, caused by the the squared excess bias in the final bin.}
\label{fi:StdevGivenMax}
\end{figure}

Figure \ref{fi:MeanArrGivenMax} plots the time averaged value of $f(t,h)$ versus $h$. The curve looks concave, smooth and possibly nearly linear for larger
values of $h$. Figure \ref{fi:StdevGivenMax} plots the square root of the time average of $Var[B(t|h)]$.
Even after time averaging, the curve is noisy in its $h$ dependence.
We note that the increase in $\int_0^1 var(t,h)dt$ may contain substantial bias for the largest point(s) in $h$. 

For each of the ten values of the high, we display both the simulation curve and the analytic curve from numerically computing the moments of \eqref{Prbxth}.
The simulated curves have the symbols overstruck on them. The point is the match of simulation with \eqref{Prbxth} is very good.
\begin{figure}[htbp]
\centering
\includegraphics[scale=.4]{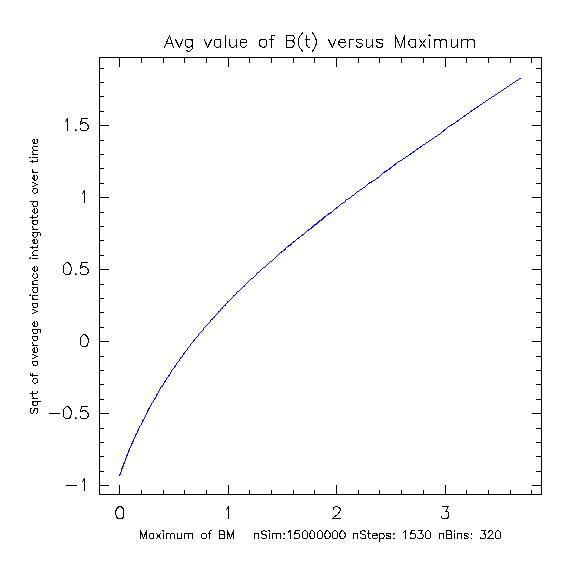}
\caption{Time average of  $E[B(t)|h]$ versus its maximum, $h$}
\label{fi:MeanGivenMax}
\end{figure}

\subsection{Feller Range \label{FellerSect}}
\begin{figure}[htbp]
\centering
\includegraphics[scale=\FSCL]{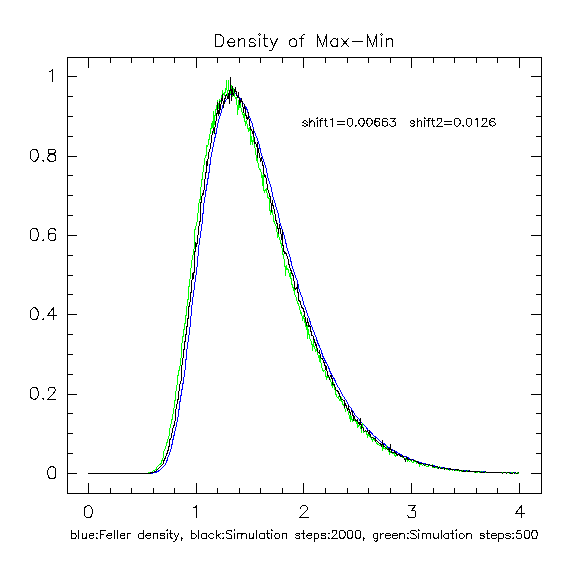}
\caption{Density of the Range: $\max_{0\le s \le 1}B(s) -\min_{0\le s \le 1} B(s)$. Computation of Feller's formula uses 400 terms.
Empirical distribution uses 4,000,000 simulations with 500 and 2000 time steps.}
\label{fi:Frangestp2000}
\end{figure}

To look at convergence, we examine the distribution of the range as a function of the number of steps
in the Brownian motion simulation. The theoretical distribution was calculated by Feller in \cite{Feller}:
The range $R_t\equiv \max_{0\le s \le t}B(s) -\min_{0\le s \le t} B(s) $
at time $t$ is distributed like $\sqrt{t}R_1$ and
the density of $R_1$ is the function 
$f(x)=8 \sum_{k=1}^{\infty} (-1)^{k+1}k^2 \phi(kx)$
where $\phi$ denotes the standard normal density and defined on $(0,\infty)$ \cite{Feller}.
As noted by Feller: "In this form it is not even obvious that the function is positive".
We compute Feller's formula for the density of the range of a Brownian motion. It converges very slowly near zero.
To evaluate $f(x=.005)$, we need between $300$ and $400$ terms. The formula is useful to compare
our Brownian motion computations with the theoretical results.
Although Feller's article is almost seventy years old, we are unaware of any previous numerical study of
its convergence or even a computation of it.
Figure \ref{fi:Frangestp2000} compares the empirical density with Feller's result. The blue curve is computed from
Feller's expansion, the black curve is the empirical density from four million realizations with 2000 time steps.
The green curve uses only 500 time steps.
We see very good agreement. The main difference is that the empirical distribution is shifted slightly to the left.
There is less than $0.1\%$ of the distribution below $range<0.7$. In Section \ref{HLCcalcSect},
the density given high and low bounds involves
an expansion in $\sum_k \exp(- k^2(h-l)^2)$. This expansion converges very quickly for vast majority of the ensemble of Brownian
paths. 

We see that the shift of the empirical distribution decreases as the step size decreases.
For a step size of .0005, the shift of the center of mass of the distribution is .0066 from the theoretical result. Using
a timestmp four times larger doubles the shift.

The distribution of the range is very small for $range<.5$ and this region is poorly approximated by the Feller expansion.
The is the clear opportunity for an asymptotic expansion in the region of small range.

\section{Figures Conditional on Close, High \label{TDCH}}

In this section, we plot the $E[B(t|c,h)]$ and $Var[B(t|c,h)]$ for a variety of different values of $(c,h)$.
Specifically, we choose quantiles (.2,.5, .8) of the bin values for the close. For our robustified grid, this corresponds
to $close = -1.011, 0.0152, 1.055$. 
In each plot, we plot the expectation $E[B(t|c,h)]$ for ten values of $h$. The ten values of $h$ are chosen to be equi-spaced in the bin coordinate
from the second bin to the second to the last bin.
We then repeat for $Var[B(t|c,h)]$. We conclude with plots for the time average of $E[B(t|c,h)]$ and $Var[B(t|c,h)]$.

For these plots, we use $1530$ time steps on each simulation for a total of 18 million simulations with 100 bins in each parameter direction. The curves overstruck by symbols are the simulation curves. The analytic formula curves have the same color but no symbol.

\subsection{ Time Dependent Mean Given Close, High}
\begin{figure}[htbp]
\centering
\includegraphics[scale=.5]{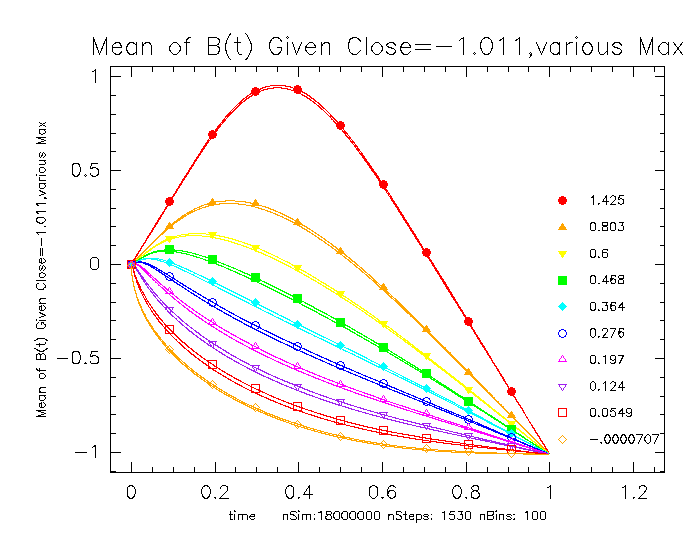}
\caption{$E[B(t |close=-1.011, various\ high)]$ where the values of the maximum are given in the legend }
\label{fi:CHLMeanGiven2Close=-1.011CHL}
\end{figure}

\begin{figure}[htbp]
\centering
\includegraphics[scale=.5]{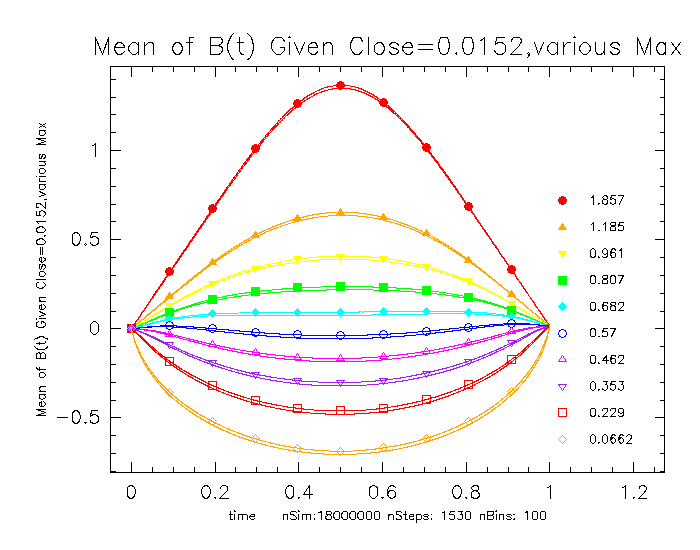}
\caption{$E[B(t |close=0.0152, various\ high)]$. The smooth curves with no symbol are given by \eqref{CHMom1} while the noisy curves
are our simulation.}
\label{fi:CHLMeanGiven2Close=0.0152CHL}
\end{figure}

Figure \ref{fi:CHLMeanGiven2Close=-1.011CHL} shows that the expectation is nearly monotonically decreasing for strongly negative values of the close
and near zero values of the high. We say nearly decreasing because we have not examined the behavior near $time=0$.
For large values of the high, the high peaks near the middle of the time interval.
Figure \ref{fi:CHLMeanGiven2Close=0.0152CHL} shows the expectation is nearly symmetric in time when the close is near zero.
\begin{figure}[htbp]
\centering
\includegraphics[scale=.5]{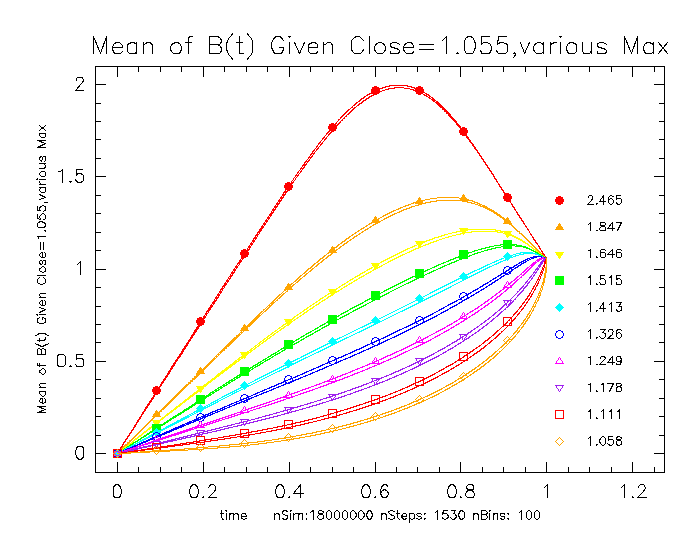}
\caption{$E[B(t |close=1.055, various\ high)]$. The values of the high are given in the legend.}
\label{fi:CHLMeanGiven2Close=1.055CHL}
\end{figure}
Figure \ref{fi:CHLMeanGiven2Close=-1.011CHL}  and Figure \ref{fi:CHLMeanGiven2Close=1.055CHL} display the following reflection
symmetry: $E[B(t |-c, h)] = E[B(1-t |c, h+c)] -c$ where $c>0$. 

\subsection{ Time Dependent Variance Given Close, High}
\begin{figure}[htbp]
\centering
\includegraphics[scale=.5]{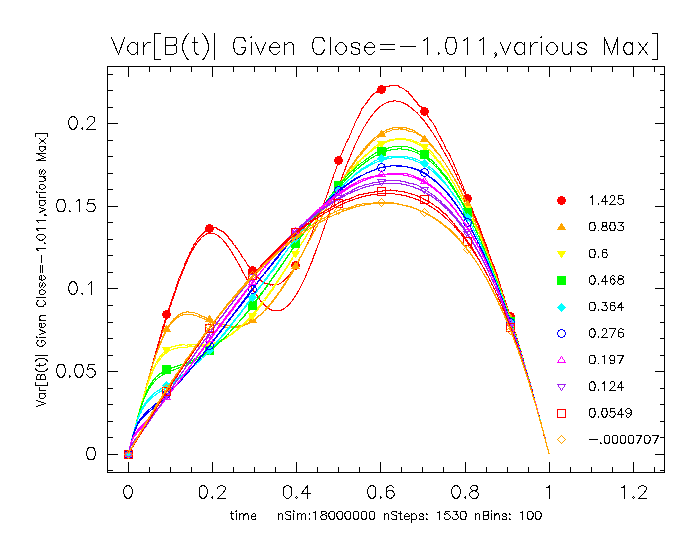}
\caption{$Var[B(t |close=-1.011, high)]$.}
\label{fi:CHLVarGiven2Close=-1.011CHL}
\end{figure}

\begin{figure}[htbp]
\centering
\includegraphics[scale=.5]{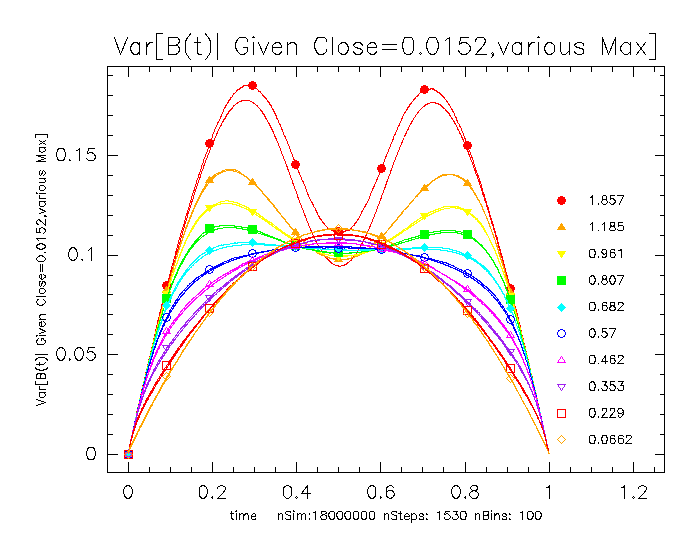}
\caption{$Var[B(t |close=0.0152, high)]$}
\label{fi:CHLVarGiven2Close=0.0152CHL}
\end{figure}

\begin{figure}[htbp]
\centering
\includegraphics[scale=.5]{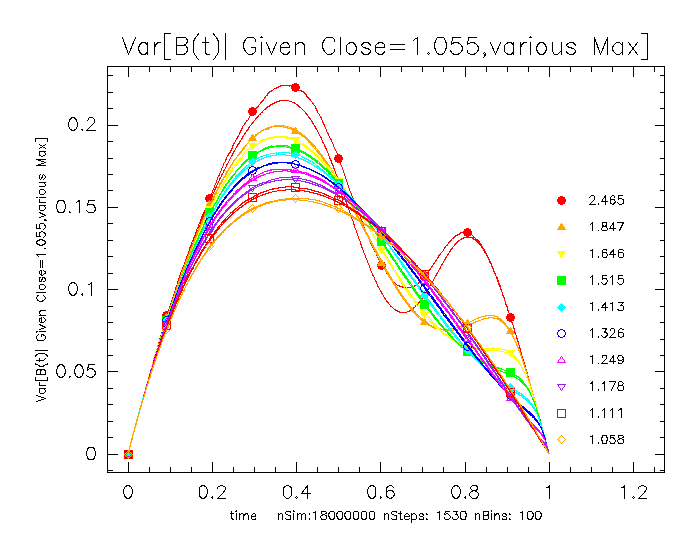}
\caption{$Var[B(t |close=1.055, high)]$}
\label{fi:CHLVarGiven2Close=1.055CHL}
\end{figure}

Figures \ref{fi:CHLVarGiven2Close=-1.011CHL}-\ref{fi:CHLVarGiven2Close=1.055CHL} display $Var[B(t |close, high)]$ for
$close = -1.011, 0.0152, 1.055$. The smooth curves with no symbol are the analytic results from \eqref{CHMom1} and \eqref{CHMom2}.
In many cases, the variance is multimodal in time.
Figure \ref{fi:CHLVarGiven2Close=-1.011CHL}  and Figure \ref{fi:CHLVarGiven2Close=1.055CHL} display the following reflection
symmetry: $Var[B(t |-c, h)] = Var[B(1-t |c, h+c)]$ where $c>0$.

\clearpage
\newpage
\section{Comparison of Theory and Simulation Given Close and High} \label{CompCHSect}

\begin{figure}[htbp]
\centering
\includegraphics[scale=.5]{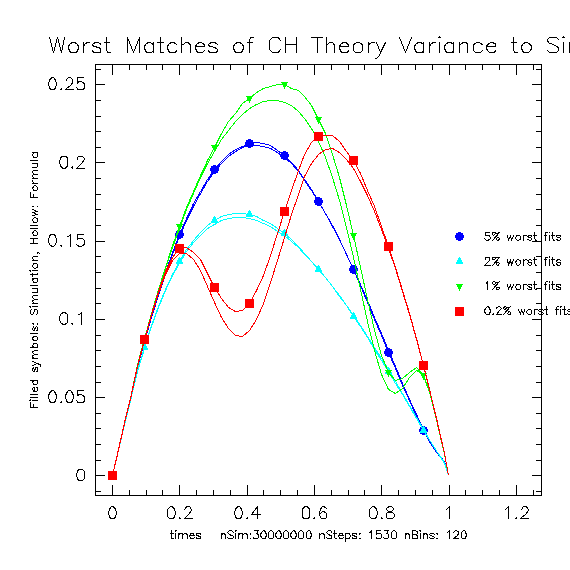}
\caption{Comparison of simulation with \eqref{CHMom2A} for four values of $(c,h)$. \quad
$5\%$ worst MSE Var:.0000000363 at close:1.996, high:2.167 \quad
$2\%$ worst MSE Var:.0000000567 at close:1.007, high:1.148 \quad
$1\%$ worst MSE Var:.00000236 at close:3.125, high:4.036  \quad
$0.2\%$ worst MSE Var:.00000277 at close:-0.836, high:1.518}
\label{fi:bo_CHVarCombo}
\end{figure}

In this section, we plot the simulation and theoretical calculation given by \eqref{CHMom1} and \eqref{CHMom2A}.
for this comparison, we use 30 million realizations each with 1530 steps. The results are then binned
in 120 bins in each direction for a total of 1.73 million bins. We compute the MSE for each bin
and sort them. We then display the fits for the worst .05, .02, .01 and .002 quantiles of the bins.
To put the curves to scale, we plot all the curves together. The curves overstruck by symbols are the simulation curves. The analytic formula curves have the same color but no symbol.
\begin{figure}[htbp]
\centering
\includegraphics[scale=.5]{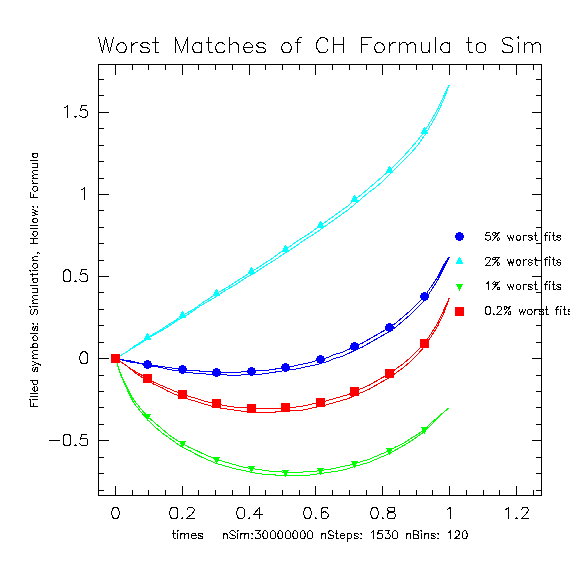}
\caption{Comparison of simulation with \eqref{CHMom1} for four values of $(c,h)$. 
$5\%$ worst MSE Mean:.0000117 at close:0.622, high:0.718 \quad
$2\%$ worst MSE Mean:.0000124 at close:1.67, high:1.739 \quad
$1\%$ worst MSE Mean:.000013 at close:-0.294, high:0.0581 \quad
$0.2\%$ worst MSE Mean:.0000136 at close:0.373, high:0.446}
\label{fi:bo_CHMeanCombo}
\end{figure}
We now display the comparisons for each bin separately. This rescales the y-axis and makes the comparison 
look worse. The differences are due to a) averaging realizations for different values of $(c,h)$;
b) discretization errors from the finite time step of the Brownian motion.
The black curve is the analytic expression while the blue curve is the the ensemble average of the simulation within
the given bin.
Figure \ref{fi:bo_CHVarCombo} compares the simulated variance in four separate bins with the
analytic expression in \eqref{CHMom2A}. Here again, we compute the squared error for each of the one million
bins. We then plot the fits for the worst .05, .02, .01 and .002 quantiles of the bins.
The worst fits for the variance have different parameters than the parameters for the worst fits to the empirical mean.
To put the curves to scale, we plot all the curves together.

\section{Figures Conditional on Close, High, Low \label{TDCHL}}
In this section, we plot the $E[B(t|c,h,\ell)]$ and $Var[B(t|c,h,\ell)]$ for a variety of different values of $(c,h,\ell)$.
Specifically, we choose quantiles (.2,.5, .8) of the bin values for the close. For our robustified grid, this corresponds
to $close = -1.011, 0.0152, 1.055$. For each value of the close, we choose three values for the high corresponding
to the $(.2,.5,.8)$ quantiles of the roubstified grid in $h$. This gives nine plots for $E[B(t|c,h,\ell)]$. 
In each plot, we plot the expectation $E[B(t|c,h,\ell)]$ for ten values of $\ell$.
We then repeat for $Var[B(t|c,h,\ell)]$.
For these plots, we use $1530$ time steps on each simulation for a total of 18 million simulations with 100 bins in each parameter direction. The curves overstruck by symbols are the simulation curves. The analytic formula curves have the same color but no symbol.

\subsection{ Time Dependent Mean Given Close, High, Low}

\begin{figure}[htbp]
\centering
\includegraphics[scale=.5]{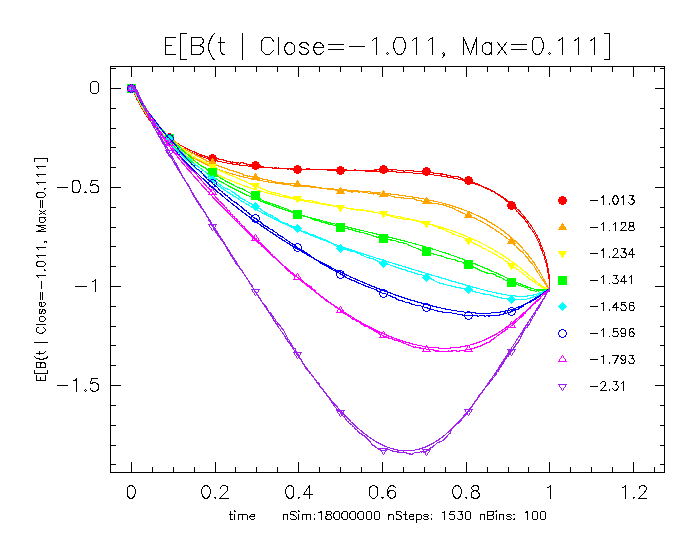}
\caption{$E[B(t |close=-1.011, high=.111,various\ low)]$ where the values of the low are given in the legend. Smaller values of the low occur on average earlier in time.}
\label{fi:CHLMeanGivenClose-1.011Max0.111CHL}
\end{figure}

\begin{figure}[htbp]
\centering
\includegraphics[scale=.5]{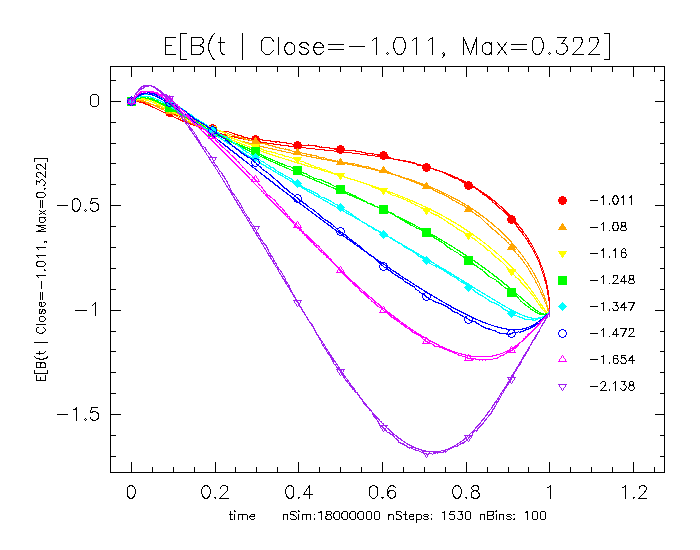}
\caption{$E[B(t |close=-1.011, high=0.322, various\ low)]$}
\label{fi:CHLMeanGivenClose-1.011Max0.322CHL}
\end{figure}

\begin{figure}[htbp]
\centering
\includegraphics[scale=.5]{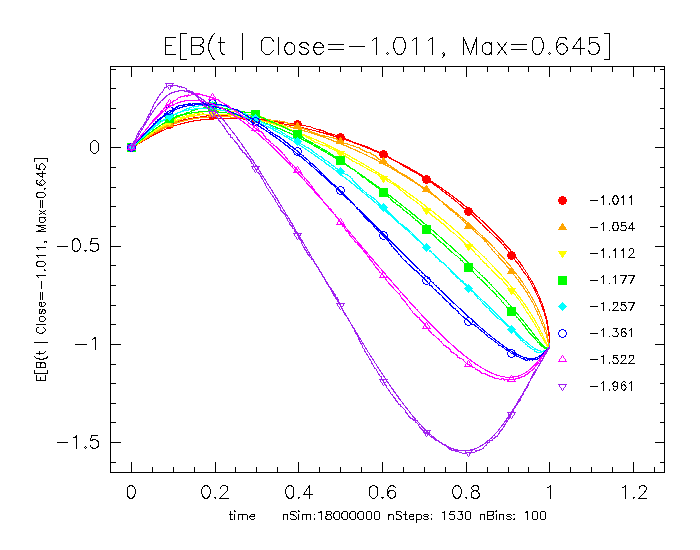}
\caption{$E[B(t |close=-1.011, high=0.645, various\ low)]$. The values of the minimum are given in the legend.}
\label{fi:CHLMeanGivenClose-1.011Max0.645CHL}
\end{figure}

\begin{figure}[htbp]
\centering
\includegraphics[scale=.5]{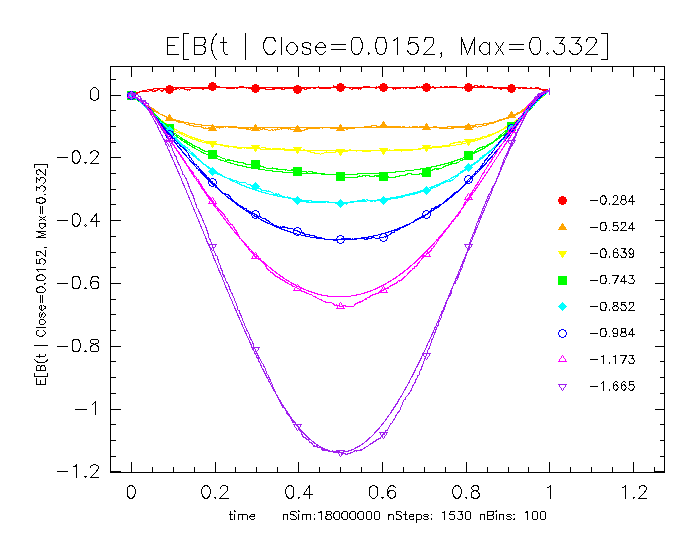}
\caption{$E[B(t |close=0.0152, high=0.332, various\ low)]$ }
\label{fi:CHLMeanGivenClose0.0152Max0.332CHL}
\end{figure}

\begin{figure}[htbp]
\centering
\includegraphics[scale=.5]{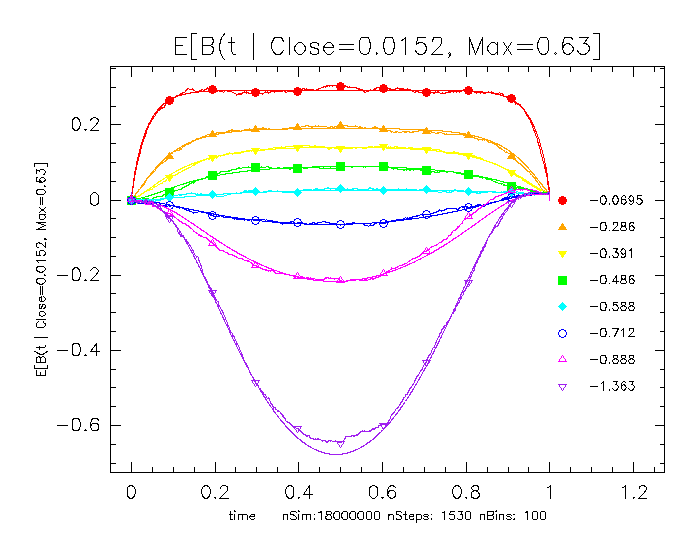}
\caption{$E[B(t |close=0.0152, high=0.63, various\ low)]$ }
\label{fi:CHLMeanGivenClose0.0152Max0.63CHL}
\end{figure}

\begin{figure}[htbp]
\centering
\includegraphics[scale=.5]{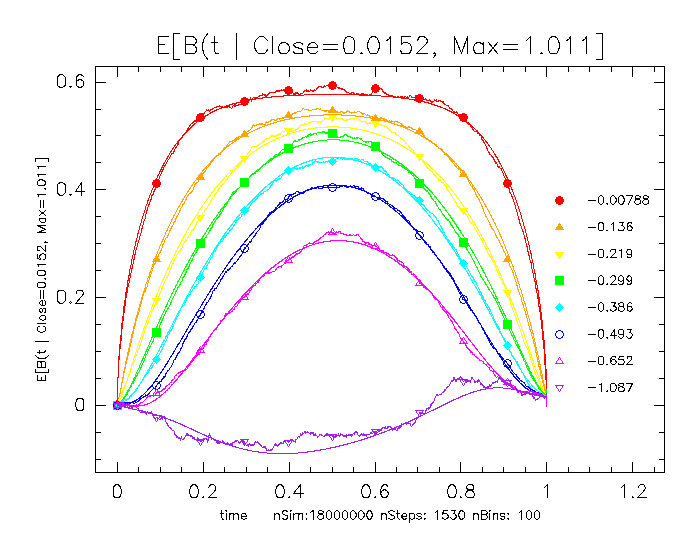}
\caption{$E[B(t |close=0.0152, high=1.011, various\ low)]$ }
\label{fi:CHLMeanGivenClose0.0152Max1.011CHL}
\end{figure}

Figures \ref{fi:CHLMeanGivenClose-1.011Max0.111CHL}- \ref{fi:CHLMeanGivenClose-1.011Max0.645CHL}
show $E[B(t|c,h,\ell)]$ for $c=-1.011$. Note that maximum of the expectation is less the expectation of the maximum.
The curves on the three plots have a similar shape as the value of the low is varied. This may indicate a somewhat weaker
dependence on high than on the low when the close equals -1. However, a stronger factor is that the curves in 'low' the low coordinate
vary more since we sample 10 values from the second smallest bin value of $\ell$ to the second largest value of $\ell$ given $(c,h)$.
\clearpage
Figures \ref{fi:CHLMeanGivenClose0.0152Max0.332CHL}-\ref{fi:CHLMeanGivenClose0.0152Max1.011CHL}
show $close=0.0152$. In this case (close near zero), the expectation is roughly symmetic.
In Figure \ref{fi:CHLMeanGivenClose0.0152Max1.011CHL}, the curves for large high and small low
are not very symmetric, but this may be due to fewer curves in the bin due to our adaptive binning.
We have also plotted  $E[B(t |c, h, \ell)]$ for $c \approx 1$. These plots exhibit the same reflection symmetry that
Figure \ref{fi:CHLMeanGiven2Close=-1.011CHL}  and Figure \ref{fi:CHLMeanGiven2Close=1.055CHL} do.
Specifically,  the reflection symmetry: $E[B(t |-c, h, \ell)] =E[B(1-t |c, h+c, \ell+c)]-c $ where $c>0$.

\newpage
\subsection{ Time Dependent Variance Given Close, High, Low}
\begin{figure}[htbp]
\centering
\includegraphics[scale=.5]{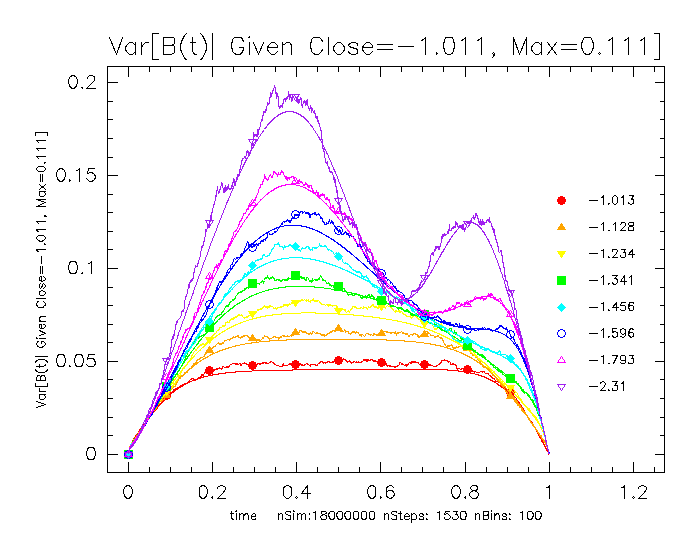}
\caption{$Var[B(t |close=-1.011, high=0.111, low)]$. The curves are roughly symmetric in time for small $|low|$
and multimodal for large $|low|$.}
\label{fi:CHLVarGivenClose-1.011Max0.111CHL}
\end{figure}

\begin{figure}[htbp]
\centering
\includegraphics[scale=.5]{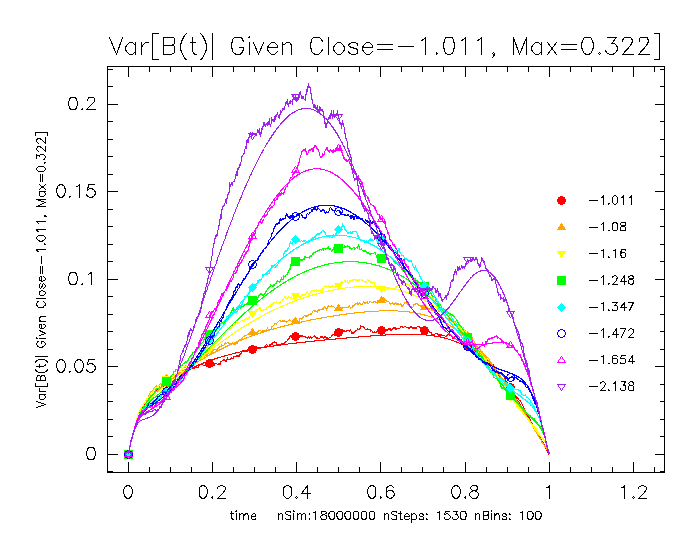}
\caption{$Var[B(t |close=-1.011, high=0.322, low)]$}
\label{fi:CHLVarGivenClose-1.011Max0.322CHL}
\end{figure}

\begin{figure}[htbp]
\centering
\includegraphics[scale=.5]{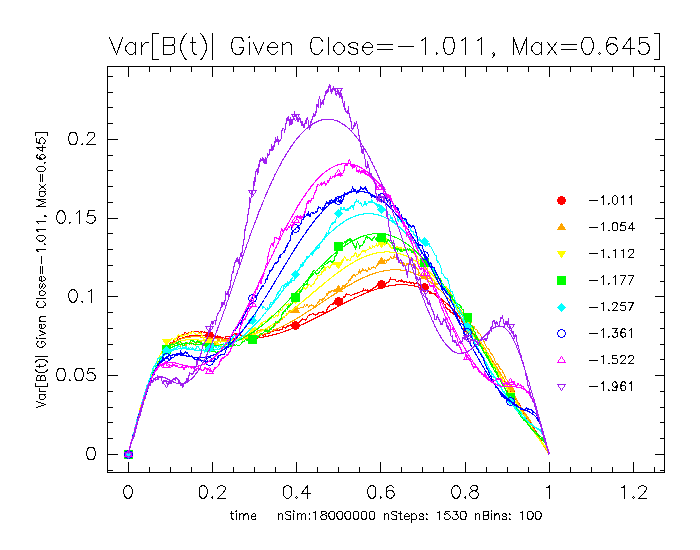}
\caption{$Var[B(t |close=-1.011, high=:0.645, low)]$}
\label{fi:CHLVarGivenClose-1.011Max0.645CHL}
\end{figure}

\begin{figure}[htbp]
\centering
\includegraphics[scale=.5]{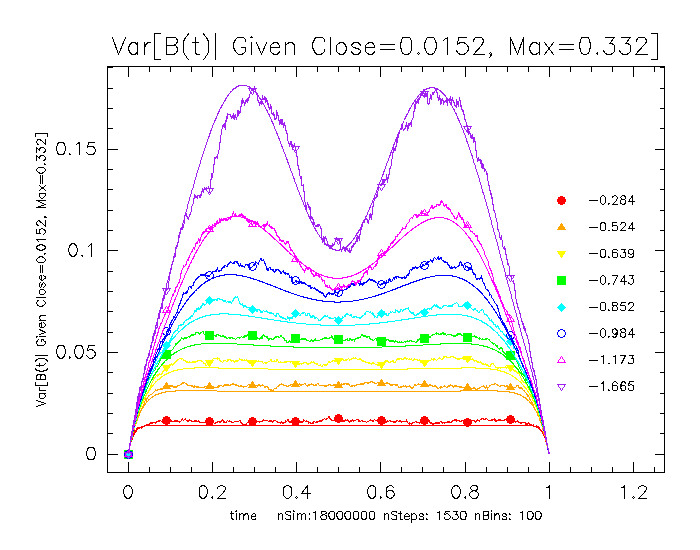}
\caption{$Var[B(t |close=0.0152, high=:0.332, low)]$:VarGivenClose0.0152 Max0.332CHL}
\label{fi:CHLVarGivenClose0.0152Max0.332CHL}
\end{figure}

\begin{figure}[htbp]
\centering
\includegraphics[scale=.5]{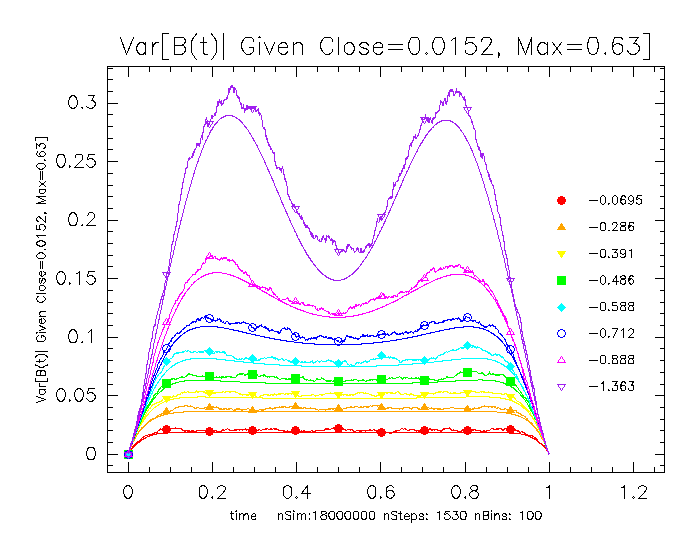}
\caption{$Var[B(t |close=0.0152, high=0.63, low)]$}
\label{fi:CHLVarGivenClose0.0152Max0.63CHL}
\end{figure}

\begin{figure}[htbp]
\centering
\includegraphics[scale=.5]{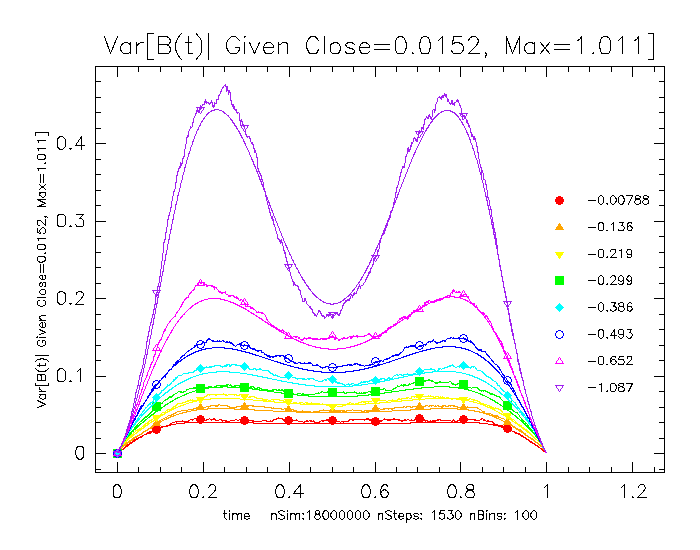}
\caption{$Var[B(t |close=0.0152, high=1.011, low)]$}
\label{fi:CHLVarGivenClose0.0152Max1.011CHL}
\end{figure}

Figures \ref{fi:CHLVarGivenClose-1.011Max0.111CHL}-\ref{fi:CHLVarGivenClose0.0152Max1.011CHL}
plot the $Var[B(t|c,h,\ell)]$ for $close = -1.011, 0.0152, 1.055$. For each value of the close, we choose
we choose quantiles (.2,.5, .8) of the bin values for the high. In many cases, the variance is multimodal in time.
The curves are much noisier because the 18 million realizations are now put into 10,000 bins instead of 100 bins.

Figures \ref{fi:CHLVarGivenClose0.0152Max0.332CHL}-\ref{fi:CHLVarGivenClose0.0152Max1.011CHL} display the empirical variance
when the close is near zero. In many cases, the variance is strongly bimodal with maxima near $t=.25$ and $t=.75$.
The $y$-axis is self-scaled. The largest uncertainties occur for $c$ near zero.
Similarly, $Var[B(t |c, h, \ell)]$ for $c \approx 1$. These plots exhibit the same reflection symmetry that 
Figure \ref{fi:CHLVarGiven2Close=-1.011CHL}  and Figure \ref{fi:CHLVarGiven2Close=1.055CHL} do/
Specifically,  the reflection symmetry: $Var[B(t |-c, h, \ell)] =Var[B(1-t |c, h+c, \ell+c)]$ where $c>0$.


\clearpage

\section{Comparison of Theory and Simulation Given Close, High and Low \label{TDCHLcomp}  \label{CompCHLSect} }

In this section, we plot the simulation and theoretical calculation given by \eqref{HLCmom0}.
For this comparison, we use 30 million realizations each with 1530 steps. The results are then binned
in 120 bins in each direction, thus a total of 1.73 million bins. We compute the MSE for each bin
and sort them. We then display the fits for the worst .05, .02, .01 and .002 of the bins.
To put the curves to scale, we plot all the curves together. 
\begin{figure}[htbp]
\centering
\includegraphics[scale=.5]{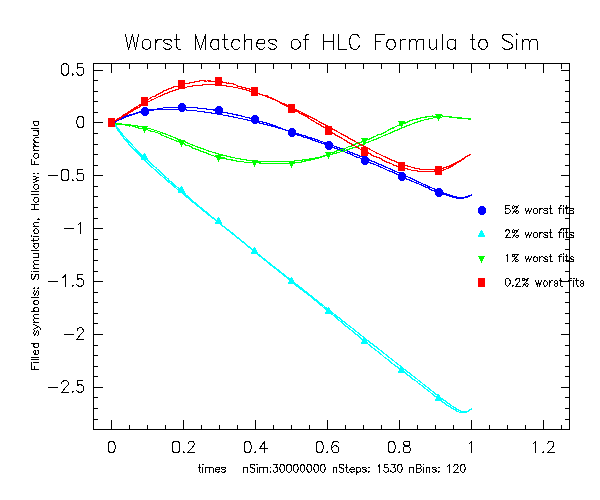}
\caption{Comparison of simulation with \eqref{HLCmom0} for four values of $(c,h,\ell)$. \quad
$5\%$ worst MSE Mean:.0000114 at close:-1.289, high:0.109, low:-1.502; \quad
$2\%$ worst MSE Mean:.000014 at close:0.836, high:0.932, low:-0.487; \quad  
$1\%$ worst MSE Mean:.0000165 at close:0.972, high:1.256, low:-1.198; \quad
$0.2\%$ worst MSE Mean:.0000257 at close:0.242, high:0.875, low:-1.056.
}
\label{fi:bo_HLCMeanCombo}
\end{figure}

\begin{figure}[htbp]
\centering
\includegraphics[scale=.5]{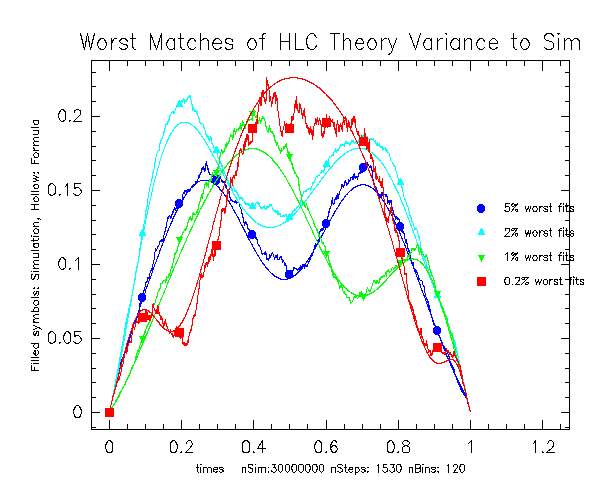}
\caption{Comparison of simulation with \eqref{HLCmom0} for four values of $(c,h,\ell)$.
$5\%$ worst MSE Var:.00000124 at close:-0.038, high:0.756, low:-0.611;
$2\%$ worst MSE Var:.0000021 at close:3.125, high:3.286, low:-0.276;
$1\%$ worst MSE Var:.00000311 at close:-3.125, high:0.307, low:-3.642;
$0.2\%$ worst MSE Var:.00000695 at close:-3.125, high:0.911, low:-3.211.
}
\label{fi:bo_HLCVarCombo}
\end{figure}

Figure \ref{fi:bo_HLCVarCombo} compares the simulated variance in four separate bins with the
analytic expression in \eqref{CHMom2A}. Here again, we compute the squared error for each of the 1.73 million
bins. We then plot the fits for the worst .05, .02, .01 and .002 quantiles of the bins.
The worst fits for the variance have different parameters than the parameters for the worst fits to the empirical mean.
To put the curves to scale, we plot all the curves together.

\section{Estimation of SP500 Prices  \label{SP500Sect}}

We now estimate the SP500 index future, ES, given only prices at the open, high, low and close.
Our real applications use the formulas in Section \ref{HLCcalcSect}  to estimate the time evolution
of series with only open, high, low and close data. We choose the SP500 because we have the
time history and can test the performance of various estimators. For K days between 2005 and 2015
we compute ten second bars between 9:30 am and 16:00 pm for a total of $2340=6*390$ prices per day.
We exclude half days. For each day, we define
$X(t,day)\equiv \log(price)(t,day) - \log(price)(t=0,day)$.
Since $X(t=0,day)=0$, we do not use the first point each day.
The volatility varies throughout the day, being larger near the beginning and end of the day.

We estimate the time dependence as 
\BEQ
\sggh(t) = \frac{1}{N_d}\sum_{day} [X(t,day) - X(t-h,day)]^2
\NEQ
where $N_d$ is the number of days in the sum. 
The time dependent volatility is independent of day. 
The sum of the volatilities,  $\sum_{j=1}^N \sggh(t_j) =.007216^2$, coresponding to
a daily volatility of $1.6$.

We define volatility time $\tau_i$ by
\BEQ \label{voltime}
\tau_k =\tau(t_k) = \sum_{i=1}^k \sggh(t_i) / \sum_{j=1}^N \sggh(t_j) 
\NEQ
In volatility time, diffusion rate of $X(\tau(t))$ is time independent and matches 
the assumptions of Brownian motion.

We score our various estimates, $\Xh(t,day)$, with the MSE in volatility time:
\BEQ \label{volscore}
MSE=  \frac{1}{N*N_d} \sum_{i,day} \sggh(t_i) [X(t_i,day) - \Xh(t_i,day)]^2
\NEQ
To normalize the MSE, we use the relative mean square error:
\BEQ \label{Rvolscore}
RMSE= \sum_{i,day} \sggh(t_i) [X(t_i,day) - \Xh(t_i,day)]^2 /  \sum_{i,day} \sggh(t_i) X(t_i,day)^2
\NEQ
We also give the mean relative squared error:
\BEQ \label{Relvolscore}
MRSE= \sum_{day}  \frac{\sum_i\sggh(t_i) [X(t_i,day) - \Xh(t_i,day)]^2}{ \sum_{i} \sggh(t_i) X(t_i,day)^2}
\NEQ
We include the estimated variance, $\sggh(t_i)$, in the loss measure because it corresponds to a time integral in volatility time.

The estimators from Section \ref{HLCcalcSect} require $\sgg$ and this must be estimated.
The simplest is the date independent estimate $\sggh_{const} = Mean[c(day)^2]$.
The Garmen Klass estimate is $\sggh_{GK}(h,\ell,c)= K_1(h-\ell)^2  -K_2[c(h+\ell)-2h\ell] - K_3 c^2 $ where $K_1=.511$,
$K_2=.019$ and $K_3=.383$ \cite{GarmenKlass, RogersSatchell}. The maximum likelihood estimator, $\sggh_{ML}(h,\ell,c)$, is based on \eqref{dChoiRoh}.
The Meillijson estimator, $\sggh_{M}(h,\ell,c)$ is given in \cite{Meillijson}. We reject the Rogers-Satchell estimator because it gives $\sggh_{R}(h=0,\ell,c)=0$
when the low happens on the close ($\ell=c$) and $h=0$. This case does occur in financial data even though it never occurs in Brownian motion.

Table 1 compares the MSE of the estimates.
Our first estimates is the are the Brownian bridge, $\Xh(t_i,day) = c_{day} t_i$, and
the second row is the Brownian volatility bridge, $\Xh(t_i,day) = c_{day} tau_i$.
The remainder of our estimate are $E[B(t|h,\ell,c)]$  as given by Theorem \ref{momThmHLC} with various plu in estimates of $\sigma$.
In Table 1, $E[B(\tau|h,\ell,c)]$ denotes using volatility time. We see the use of volatility time only slightly improves the fit. For these fits, we estimate the volatility every day separately, but use volatility time calculated for
the whole data set.
The maximum likelihood estimate, $\sggh_{ML}(h,\ell,c)$, does slightly better than $\sggh_{GK}(h,\ell,c)$ with $\sggh_{M}(h,\ell,c)$ coming in third.
Using the information from $(h,\ell,c)$ improves reduces the error to  $56\%$ of the error of the Brownian bridge.

\begin{table}[ht] \label{tab:EstPerfTab}
\caption{Performance of Various Estimators on SP500 Data }
\centering
\begin{tabular}{cccc}
\hline 
Estimator& Sigma Estimate & RMSE & MRSE \\
\hline
Bridge& N.A. & 0.37089 & 0.62048\\
Bridge Vol Time &  N.A.& 0.36332 & 0.63496 \\
Thm \ref{momThmHLC} $E[B(t|h,\ell,c)]$    & $\sigma_{GK}$ & 0.21210 & 0.40169 \\
Thm \ref{momThmHLC} $E[B(\tau|h,\ell,c)]$ & $\sigma_{GK}$ & 0.20275 & 0.38737\\
Thm \ref{momThmHLC} $E[B(t|h,\ell,c)]$    & $\sigma_{ML}$ & 0.20928 & 0.39615\\
Thm \ref{momThmHLC} $E[B(\tau|h,\ell,c)]$ & $\sigma_{ML}$ & 0.20368 & 0.38517\\
Thm \ref{momThmHLC} $E[B(\tau|h,\ell,c)]$ & $\sigma_{M}$ &  0.23754 & 0.4084\\
\hline \hline
\end{tabular}
\end{table}

\section{Distribution Given High and Low \label{HLSect} }
We evaluate  the distribution $p(x, t,h,\ell)$ by integrating over the closing value in $p(x,t,h,\ell,c)$ using \eqref{CHLSum1}. As before, the limits of integration, $H$ and $L$, are to be set to $h$ and $\ell$ after differentiation.

\begin{theorem}
  Let $G_{HL}(x,t, h, \ell) \equiv P(B(t)=x, \ell \leq B(s) \leq h | {\rm for}\ s \in [0,1] )$
\BEQ \label{HLGen}
G_{HL}(x,t, h,\ell) \equiv \int_L^H G(x,t, h,\ell, c) dc = Q(x,t, h, \ell) \sum_{k>-\infty}^{\infty}\left[R_{1k} (x,t, h, \ell;H,L) -R_{2k}(x,t, h, \ell)\right]
\NEQ
where $Q(x,t, h, \ell)$ is defined in \eqref{ChoiRoh} and
\begin{align*}\label{HLGen1}
R_{1k} &= \frac{1}{2}\left[erf(\frac{H-x -k \Delta}{\sqrt{2(1-t)}\sigma}) - erf(\frac{{L-x -k \Delta} }{\sqrt{2(1-t)}\sigma})\right] , \  \\
R_{2k} &=\frac{1}{2} \left[erf(\frac{H+x -2h +k \Delta}{\sqrt{2(1-t)}\sigma}) - erf(\frac{{L+x -2h +k \Delta} }{\sqrt{2(1-t)}\sigma})\right] \ \ .
\end{align*}
The density satisfies $p(x;t,h,\ell) = - \lim_{H\rightarrow h, L \rightarrow \ell}\partial_{\ell}\partial_h G(x,t, h,\ell;H,L)$.
\end{theorem}

\Prf \eqref{HLGen} is $Q(x,t, h, \ell)\int_L^H Q_R(x,t, h,\ell, c) dc$, integrated term by term. We again use the articificial limits
of $L$ and $H$ to indicate that the limits should not be differentiated in evaluating the density.
\qed

To get the density conditional on the high and low, we divide $P(x;t,h,\ell)$ by $p(h,\ell)$ as given by \eqref{densHL}
The theoretical values of  $E[B(t)|h,\ell]$ and  $Var[B(t)|h,\ell]$ can be calculated by computing moments with respect to \eqref{Prbxth}.
Unfortunately, we have not found a tractable analytic form from the integrals and therefore we compute them numerically.
We display the simulation results for $E[B(t)|h,\ell]$ for a small value of $h=.304$, the median value of $h=.816$ and a large value of $h=1.572$.
\begin{figure}[htbp]
\centering
\includegraphics[scale=.5]{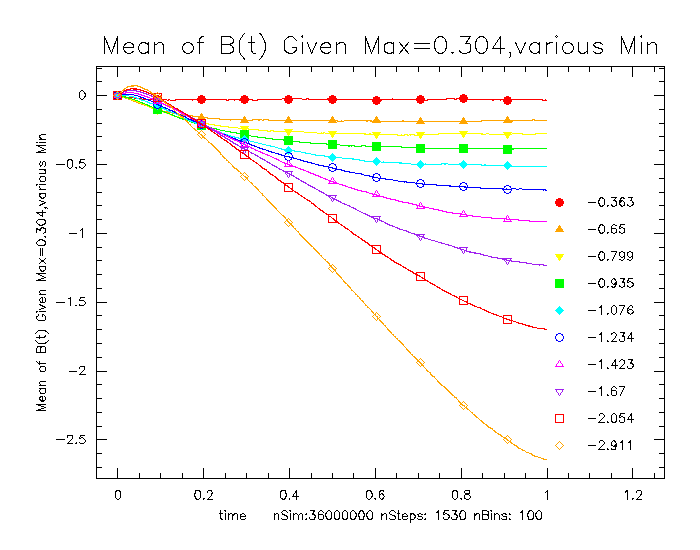}
\label{fi:HLMeanGivenMax.304}
\end{figure}
\begin{figure}[htbp]
\centering
\includegraphics[scale=.5]{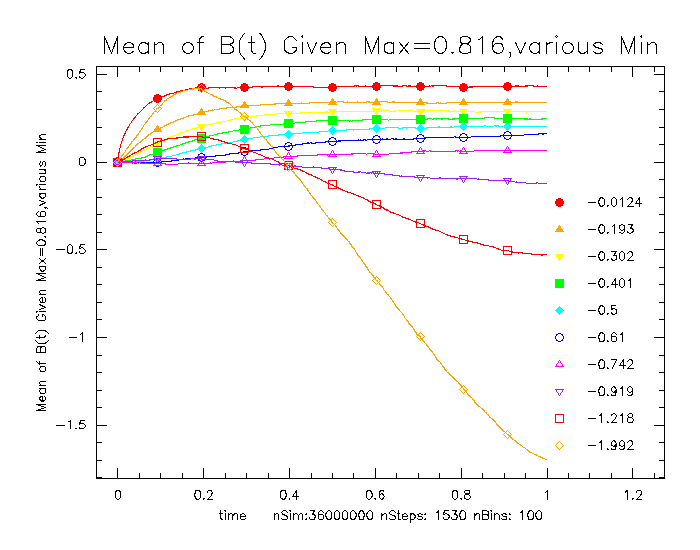}  
\caption{$E[B(t |high=0.816, various\ low)]$. If the $h> |\ell|$, the maximum occurs after the minimum.
If the $h< |\ell|$, the minimum occurs first.}
\label{fi:HLMeanGivenMax.816}
\end{figure}
\begin{figure}[htbp]
\centering
\includegraphics[scale=.5]{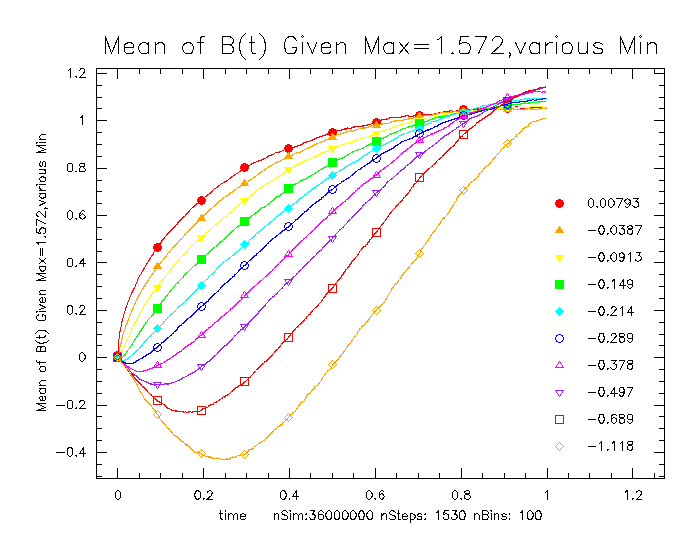}  
\caption{$E[B(t |high=1.572, various\ low)]$}
\label{fi:HLMeanGivenMax1.572}
\end{figure}

\section{Summary \label{Conclude}}

By calculating $E[B(t)|\max, \min, \close]$, we are able to interpolate in time any dataset where only the open, high, low and
close are given. In practice, we interpolate on the log scale using the logarithms of the open, high, low and close.
For most applications, we are interested in relative price chances so the log scale is appropriate. If one is truly interested in
the actual price, our formulas need to be modified for log Brownian motion.

\begin{figure}[htbp]
\centering
\includegraphics[scale=.5]{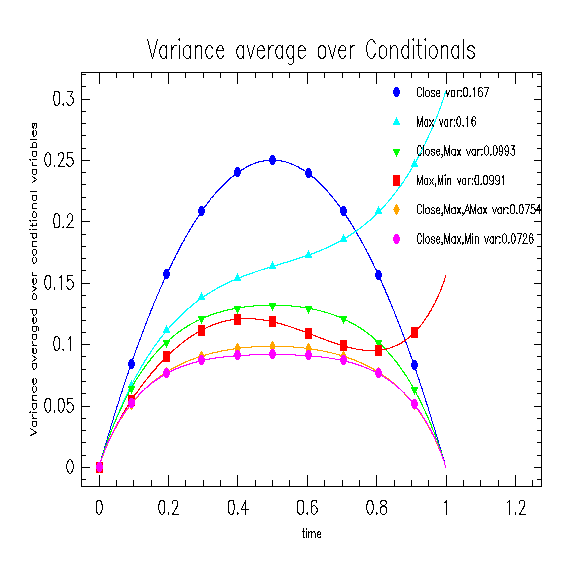}
\caption{Time Dependence of Ensemble Averaged Variance Given Conditional Variables}
\label{fi:AvgVarGivenBIG}
\end{figure}

Our simulations have calculated the ensemble average of the mean square error
in Brownian motion for a variety of different given statistics. The time dependence of the
variance is displayed in Figure \ref{fi:AvgVarGivenBIG}. In Figure \ref{fi:AvgVarGivenBIG},
\BEQ \label{Vens}
V(t|h,\ell,c)\equiv \int Var[B(t|h,\ell,c)] dp(h,\ell,c) \ ;\
V(t|h,c)\equiv \int Var[B(t|h,c)] dp(h,c) \ .
\NEQ
We ensemble average the variance expressions over all paths. For a given value of the statistics, $h$ or $(h,c)$ or $(h,\ell,c)$,
the results of the previous sections should be used for a more accurate evaluation of the variance.
The fifth curve in Figure \ref{fi:AvgVarGivenBIG} is the case when the location of the maximum is specified
in addition to $(h,c)$. This is borrowed from \cite{RiedelArgMx}.
The variance is symmetric in time when final value, $c$, is specified.
If just the high or the high and low are specified, the variance is nonmonotonic in time.

The time averaged variance in Figure \ref{fi:AvgVarGivenBIG} is presented in Table 2. 
The values for Table 2 are from the simulation. 
We plan to compute these ensemble averages using the analytic results in Sections \ref{HCcalcSect}
and \ref{HLCcalcSect}.
\begin{table}[ht] \label{tab:EnsVarTab}
\caption{Time Averaged Variance by Givens }
\centering
\begin{tabular}{ccc}
\hline 
Givens& Var& Var*6 \\
Start point only& 1/2 & 3 \\
Close & 1/6 &1\\
High &0.1602  & .9612\\
ArgMax & .2487 & 1.492\\
Close, High & 0.0990&.5938\\  
Close, ArgMax & 0.1037 &.6222\\
High, Low & 0.09911& .5947\\
Argmax, High & 0.11585 &0.6951 \\
ArgMax, ArgMin & 0.1574 & 0.9444\\
Close, High, Low & 0.0701 &.4204\\ 
Close, High, ArgMax &  0.07535&.4521\\ 
\hline \hline
\end{tabular}
\label{tabVar}
\caption{Expected time average variance reduction. We multiply the variance by $6$ in the third column to compare with knowing only the final value, $c$. Here $argMax$ is the first location of the maximum of $B(t)$. The results that contain
$argMax$ are taken from \cite{RiedelArgMx}.}
\end{table}
Table 2 shows that using the open, high, low and close reduces the variance to just $14\%$ of the variance using only
the initial value or only the final value. This shows that the use of only the open, high, low and close
in chartist forecasting \cite{Morris} keeps most of the information about the time history of the process.
Table 2 answers interesting questions like is it better to know the maximum or the final value of the
Brownian motion to predict $B(t)$ in $[0,1]$. By a ratio of $0.96$ to $1$,  
it is slightly better to know the high than the closing value.
Similarly, rows 5-8 of Table 2 show that it is better to know the close and the high than
the high or the low or the close and time of the high. The last two rows show that it is better to know the close, high and low than to know the close, high and time of the high.
Finally, we see the expected variance when using all of the open, high, low and close is just $42\%$ of the the variance from using just the open and close.

Table 1 shows the performance of our estimator on the log of the SP500 price. For the financial data, we have a MSE reduction of $54\%$ over the Brownian bridge.
This is a significant improvement but it is not as good as the theoretical value of $42\%$. The reason is clear.
For real world data, we need to estimate $\sigma$ whileour theorems and simulations have $\sigma$ given.
Also the volatility time varies from day to day in practice.

Our moment expression for $M_1(h,\ell,c)$ and  $M_2(h,\ell,c)$ in Theorem \ref{momThmHLC} and Section \ref{AppSimp} are two dimensional
sums over Gaussians terms. We are unable to collapse the two dimensional sum over $j$ and $k$ to a single infinite sum as was possible in the $m=0$ case
of Corollary \ref{AppM0}. The double expansion for $j$ and $k$ times $\Delta=(h-\ell)$  converges quickly for all but the set of Brownian paths where
$\Delta$ is very small. The Feller distribution of \ref{FellerSect} shows that the measure of the  small $\Delta$ paths is very small.

\section{Appendix: Integral Evaluations}

\subsection{Close and High \label{AppCH}}
We now evaluate the integrals $M_m$ in \eqref{CHMomi}-\eqref{CHMomiC}.  
Set $r\equiv 2h-c$,  $z_2= h-\mu_2= -z_3 = h -(2h-c)t$ and $z_4=h-\mu_4=ct-h$. 
Checking the normalization 
\BEQ
\label{CHNorm1}
M_0=\left[- \sum_{i=2}^4 s_i \tau_i f_i(x=h)\right] + 2(2h-c)\int_{-\infty}^H (f_2 +f_3) dx = p(h,c) \ .
\NEQ

For $m=1$, \eqref{CHMomiB} reduces to  
\BEQ
\label{CHMom1A}
M_1\equiv \sum_{i=2}^4 s_i \psi_i \int_{-\infty}^{h-\mu_i}\left[\tau_i -\frac{2(2h-c)}{\sigma^2}(x+\mu_i)(1-\delta_{i,4})\right] \phi_{\sigma_t^2}(x)dx
\NEQ
\BEQ
\label{CHMom1Ac}
=\sum_{i=2}^4 s_i \psi_i [\frac{\tau_i}{2} -\frac{r \mu_i}{\sigma^2} (1-\delta_{i,4})][1 +erf(\frac{h-\mu_i}{\sqrt{2}\sigma_t})] -  4 r t(1-t) \phi_{\sgg}(r) \phi_{\sgg_t}(h-rt)
\NEQ
\BEQ
\label{CHMom1Ad}
=\phi_{\sgg}(c)[1 +erf(\frac{ct-h}{\sqrt{2}\sigma_t})] + \phi_{\sgg}(r)[\frac{2hr}{\sgg} -1 +p_{h,r,t} erf(\frac{h-rt}{\sqrt{2}\sigma_t})] -  4 r t(1-t) \phi_{\sgg}(r)  \phi_{\sgg_t}(h-rt)  \ .
\NEQ
Here we define $r= r(h,c)=2h-c$ and use $\mu_2= rt$, $g_2=r^2/2$, $\mu_3=2h-rt$, $g_3=r^2/2$.
$p_{h,r,t}= (\tau_3 -\tau_2)/2 +\frac{r(\mu_2-\mu_3)}{\sgg} = (1-2t) +\frac{2r(rt-h)}{\sgg}$.
For the second moment, \eqref{CHMomiC} reduces to
\BEQ
\label{CHMom2A}
M_2\equiv \sum_{i=2}^4 s_i \psi_i \int_{-\infty}^{h-\mu_i}\left[2(x+\mu_i) \tau_i -\frac{2(2h-c)}{\sigma^2}(x+\mu_i)^2(1-\delta_{i,4})\right] \phi_{\sigma_t^2}(x)dx 
\NEQ
\BEQ
\label{CHMom2Aa}
= \sum_{i=2}^4 s_i \tau_i \psi_i
\left[\mu_i[1 +erf(\frac{h-\mu_i}{\sqrt{2}\sigma_t})] -2\sgg_t \phi_{\sgg_t}(h-\mu_i)\right]
\NEQ
\BEQ
+ \frac{2(2h-c) \phi_{\sgg}(r) }{\sgg} \left[ \sum_{i=2}^3 \frac{(\mu_i^2+\sigma_t^2)}{2} 
   \left[1 + erf(\frac{h-\mu_i}{\sqrt{2} \sigma_t})  \right] - 4h\sigma_t^2 \phi_{\sgg_t}(h-\mu_2)  \right] 
\NEQ
\BEQ
\label{CHMom2Ab}
= 2 \mu_4 \phi_{\sgg}(c)[1 +erf(\frac{ct-h}{\sqrt{2}\sigma_t})] -
\phi_{\sgg}(r)\left[ (\tau_2 \mu_2 + \tau_3 \mu_3) + (\tau_2 \mu_2 - \tau_3 \mu_3) erf(\frac{h-rt}{\sqrt{2} \sigma_t}) \right]
\NEQ
\BEQ
+ 2r \phi_{\sgg}(r) \left( \left[t(1-t) +\frac{\mu_2^2+\mu_3^2}{2\sgg} +\frac{\mu_2^2-\mu_3^2}{2\sgg} erf(\frac{z_2}{\sqrt{2} \sigma_t})\right] - 4 h t(1-t) \phi_{\sgg_t}(z_2) \right) \ =
\NEQ
\BEQ
\label{CHMom2Ac}
 2 (2h-ct)\phi_{\sgg}(c) [1 +erf(\frac{ct-h}{\sqrt{2}\sigma_t})] - \phi_{\sgg}(r) \left[q_3(h,c,t)  + q_4(h,c,t) erf(\frac{h-rt}{\sqrt{2}\sigma_t})\right]
\NEQ
\BEQ
+ 2r \phi_{\sgg}(c) \left( \left[t(1-t) + \frac{q_5(h,c,t)}{\sgg}+\frac{2h(rt-h)}{\sgg} erf(\frac{h-rt}{\sqrt{2}\sigma_t})\right] - 4 h t(1-t) \phi_{\sgg_t}(h-rt)
   \right) \ 
\NEQ
where we define $q_3(h,c,t)\equiv\mu_2 \tau_2 +\mu_3 \tau_3 =2(rt^2 +(1-t)(2h-rt) )$, $q_4(h,c,t)\equiv \mu_2 \tau_2 - \mu_3 \tau_3 =2(rt-2h(1-t)) $, $\mu_2^2+\mu_3^2=r^2t^2+ (2h-rt)^2$ and $\mu_2^2-\mu_3^2= 4h(rt-h)$. Let $q_5(h,c,t)\equiv (\mu_2^2+\mu_3^2)/2 = 2 h^2 + r^2t^2 - 2hrt = h^2 + (h-rt)^2$.


For \eqref{CHMom2Ab}, we use
\BEQ
\label{ch1}
\int_{-\infty}^{H-\mu} (x+\mu) \phi_{\sgg_t}(x)  = \frac{\mu}{2} \left[erf(\frac{H-\mu}{\sqrt{2} \sigma_t}) + 1 \right] - \sigma^2  \phi_{\sigma_t^2}(H-\mu) \ ,
\NEQ
\BEQ
\label{ch2}
\int_{-\infty}^{H-\mu}(x+\mu)^2 \phi_{\sgg_t}(x)  = \frac{\mu^2+\sigma^2}{2} \left[erf(\frac{H-\mu}{\sqrt{2} \sigma_t}) + 1 \right] - \sigma^2 (H+\mu) \phi_{\sgg_t}(H-\mu) \ .
\NEQ
These formulas have been verified by numerically integrating the moments of $\partial_h F(x,t,h,c)$ from $-\infty $ to $h$.
This completes the proof of Theorem \ref{momThm}. \qed

\subsection{$G_{mn}$ Evaluation \label{GmnEval}}

We evaluate the integrals, $G_{mn}$, of \eqref{Gmn}.
We define the scaled $erf$ function, $E_{\sigma}(x)\equiv .5* erf(x/\sqrt{2}\sigma)$. For $m=0$,
$G_{00}(\mu,h,\ell) = \left[E_{\sigma}(h-\mu)- E_{\sigma}(\ell-\mu)\right]$
and $G_{01}(\mu,h,\ell) = \sigma^2 \left[\phi_{\sigma^2}(\ell-\mu)- \phi_{\sigma^2}(h-\mu)\right]$.
\BEQ
\label{G10}
G_{10}(\mu,h,\ell,\sigma) = 
\mu \left[ E_{\sigma}(h-\mu)-E_{\sigma}(\ell-\mu) \right] +\sigma^2 \left[\phi_{\sigma^2}(\ell-\mu)- \phi_{\sigma^2}(h-\mu)\right] 
\NEQ
\BEQ
\label{G20}
G_{20}(h,\ell)   
= (\sigma^2+\mu^2) \left[ E_{\sigma}(h-\mu)-E_{\sigma}(\ell-\mu) \right] +
\sigma^2 \left[(\ell+\mu)\phi_{\sgg}(\ell-\mu)- (h+\mu)\phi_{\sgg}(h-\mu)\right] 
\NEQ
\BEQ
\label{G11}
G_{11}(\mu,h,\ell,\sigma) 
= \sigma^2 \left[ E_{\sigma}(h-\mu)-E_{\sigma}(\ell-\mu) \right] + \sigma^2 \left[\ell \phi_{\sigma^2}(\ell-\mu)- h \phi_{\sigma^2}(h-\mu)\right] \ .
\NEQ

\subsection{Centering at the Lower Limit \label{AppLowL}}

To simplify the lower limit values at $\ell$ in \eqref{HLCM1g}, we need to define the analog of $\mu_{ijk}$
except that the definitions are centered at the lower limit.
Let $\nu_{1,j,k}= ct+\vt_{j,k}\Delta$, $\gtl_1= (c-w_{j,k}\Delta)^2/2\sgg$, $\nu_2= (2\ell-c)t+v_{j,k}\Delta$, $\gtl_2=(2\ell-c-\wt_{j,k}\Delta)^2/2\sgg$, $\nu_3=2\ell(1-t) +ct-v_{j,k}\Delta$ and $\nu_4= 2\ell-ct-\vt_{j,k}\Delta$. Of course, $\gtl_4=\gtl_1= g_1$ and $\gtl_3=\gtl_2$.
The analog of \eqref{CHLSum2}, centered at the lower limits, is 
\BEQ \label{CHLSum2L}
G(x,t, h,\ell, c)=\sum_{j,k>-\infty}^{\infty} \sum_{i=1}^4 \frac{s_{i}}{\sqrt{2\pi}\sigma}\phi_{\sigma_t^2}(x- \nu_{ijk}) e^{-\gtl_{ijk}(h,\ell)}, \
\NEQ
We further define
$\taut_{ijk} \equiv \partial_{h} \nu_{ijk}$ and $\tauth_{ijk} \equiv \partial_{\ell} \nu_{ijk}$. Thus $\taut_{1jk} =\vt_{j,k}=-\tauth_{1jk}$, $\taut_{2jk} =v_{j,k}$, $\tauth_{2jk} =2t-v_{j,k}$, $\taut_{3jk}=-v_{j,k}$, $\tauth_{3jk} = 2(1-t)+v_{j,k}$, $\taut_{4jk} =-\vt_{j,k}$ and $\tauth_{4jk} =2+\vt_{j,k}$.  
Finally, we need
$A^{\ell}_{ijk} =  \taut_{ijk}  \tauth_{ijk}$, $B^{\ell}_{ijk} = [\taut_{ijk}\partial_{\ell} \gtl_{ijk}+\tauth_{ijk}\partial_{h} \gtl_{ijk}]$ and
$C^{\ell}_{ijk} =\Gamt_{ijk}+\sigma_t^{-2} \taut_{ijk}  \tauth_{ijk} $ and $\Gamt_{ijk}\equiv -\partial_h \gtl_{ijk}\partial_{\ell} \gtl_{ijk} + \partial_{\ell}\partial_h \gtl_{ijk}$.
Note $\Gamt_{1jk}=\Gamma_{1jk}$ and  $\Gamt_{2jk}=  (2*\gtl_{2jk} -1) \wt_{jk}(\wt_{jk} +2)/\sgg $.


\subsection{Further simplification of Theorem \ref{momThmHLC}\label{AppSimp}}

We now simplify \eqref{HLCM1g} by summing the Gaussian terms over $i$. The terms involving the error function do not simplify much
and are left as in Theorem \ref{momThmHLC}.
\begin{corollary}
Theorem \ref{momThmHLC} may be re-expressed as
\BEQ \label{HLCMg2}
M_m= \sum_{j,k} U_{jk}^{m,h} f_{1jk}(h-\mu_{1jk})  - U_{jk}^{m,\ell} f_{1jk}(\ell-\mu_{1jk}) +
\sum_{ijk} s_{i} e_{ijk}^{(m)} \psi_{ijk} R(h,\ell,\mu_{ijk})
\NEQ
where $R(h,\ell,\mu_{ijk}) = E_{\sigma_t}(h-\mu_{1jk})-E_{\sigma_t}(\ell-\mu_{1jk}).$   
The coefficients, $e_{ijk}^{(m)}$ are defined below \eqref{HLCM1g}.
The coefficients satisfy
\BEQ
U_{jk}^{1h} = \sum_i s_i A_{ijk} + 2(\Gamma_{2jk}-\Gamma_{1jk})\sgg_t =\Ab_{jk}  + 2(\Gamma_{2jk}-\Gamma_{1jk})\sgg_t 
\NEQ
where $\Ab_{jk}=\sum_i s_i A_{ijk} = [32 jk +8(j-k)]\ t(1-t)$. For the lower limit,
\BEQ
U_{jk}^{1\ell} = \sum_i s_i A^{\ell}_{ijk} + 2(\Gamt_{2jk}-\Gamma_{1jk})\sgg_t =\Ab^{\ell}_{jk} + 2(\Gamt_{2jk}-\Gamma_{1jk})\sgg_t
\NEQ
where $\Ab^{\ell}_{jk}=\sum_i s_i A^{\ell}_{ijk} = [32 jk-8(j-k)]\ t(1-t)$. For the second moment,
\BEQ
U_{jk}^{2h} = 2h\sum_i s_i A_{ijk}- 2 \sgg_t \sum_i s_i B_{ijk} + 4 h\sgg_t(\Gamma_{2jk}-\Gamma_{1jk})  = 2 h  \Ab_{jk}+ \Bb_{jk}  + 4h \sgg_t(\Gamma_{2jk}-\Gamma_{1jk})
\NEQ
\BEQ
U_{jk}^{2\ell} = 2\ell\sum_i s_i A^{\ell}_{ijk}- 2 \sgg \sum_i s_i B_{ijk}^{\ell} + 4 \ell\sgg(\Gamt_{2jk}-\Gamma_{1jk})  = 
2 \ell  \Ab^{\ell}_{jk} +\Btb_{jk}  + 4 \ell\sgg(\Gamma^{\ell}_{2jk}-\Gamma_{1jk})
\NEQ
where $\Bb_{jk}= \sum_i s_i B_{ijk}= \frac{8}{\sgp} [-4jk\Delta +c*j -h(j-k)] $
and $\Btb_{jk}= \sum_i s_i B^{\ell}_{ijk}=\frac{8}{\sgp} [4jk\Delta -c*j +\ell(j-k)] $.

\end{corollary}

\Prf 
We begin with
\BEQ \label{AbCalc}
\Ab_{jk}\equiv\sum_i s_i A_{ijk} = 2(v_{j,k}^2-\vt_{j,k}^2) + 2\vt_{j,k} +2(2t-1)v_{j,k}= [32 jk +8(j-k)]t(1-t) \  \ ,
\NEQ
\BEQ \label{AellbCalc}
\Ab^{\ell}_{jk} \equiv\sum_i s_i A^{\ell}_{ijk}= 2(v_{j,k}^2-\vt_{j,k}^2) - 2\vt_{j,k} -2(2t-1)v_{j,k}=[32 jk -8(j-k)]t(1-t) \  \ . 
\NEQ
To sum $s_i B_{ijk}$, we begin with the pairs $B_{ijk}$ yields
\begin{align} 
B_{1jk}+B_{4jk} &= (\tau_{1jk}+\tau_{4jk})\partial_{\ell} g_{1jk}+(\tauh_{1jk}+\tauh_{4jk})\partial_{h} g_{1jk}]
  = 2 \partial_{\ell} g_1   \ \ , \nonumber \\
B_{2jk}+B_{3jk} &=  (\tau_{2jk}+\tau_{3jk})\partial_{\ell} g_{2jk}+(\tauh_{2jk}+\tauh_{3jk})\partial_{h} g_{2jk}]
  = 2  \partial_{\ell} g_2 \ . \label{Bbar} 
\end{align}
Thus
\begin{align} \label{BbCalc}
\Bb_{jk}&\equiv\sum_i s_i B_{ijk} =2[\partial_{\ell} g_1 -\partial_{\ell} g_2] =
\frac{2}{\sgg}\left[(c-w_{jk} \Delta) w_{jk} - (2h-c - \wt_{jk} \Delta)\wt_{jk}\right] \\
&=\frac{2}{\sgg}\left[(\wt_{jk}^2-w_{jk}^2)\Delta + c  w_{jk} -(2h-c) \wt_{jk}\right] = \frac{-32jk \Delta -8(h-c)j+ 8hk}{\sgg}
\end{align}
Similarly,
\BEQ \label{BellbCalc}
\Btb_{jk}\equiv \sum_i s_i B^{\ell}_{ijk}=2\partial_{h} (\gtl_1 -\gtl_2) = \frac{32jk \Delta +8(\ell-c)j- 8\ell k}{\sgg}
\NEQ
\qed

It is possible to make small additional simplifications of \eqref{HLCMg2}, but the resulting moment computations are
not much simpler than \eqref{HLCMg2}. The computation remains a two dimensional infinite sum.

\subsection{High, Low, Close Boundary Terms \label{AppHLCBC} }

We begin by showing the boundary terms in \eqref{CHLMomGenM1} vanish:
$\Bc_{ijk}(x=H) = - h^m [(A_{ijk}/\sigma_t^2) \partial_{x}f_{ijk}(h)$ $+ B_{ijk}f_{ijk}(h)] $.  To show these boundary terms vanish, we note
$f_{ijk}(h) = f_{1jk}(h)$ , $\partial_{x}f_{ijk}(h) = (\mu_{ijk}-h)f_{1jk}(h)/\sgp_t $. Thus $\partial_{x}f_{4jk}(h) =-\partial_{x}f_{2jk}(h)$.
Also note that $0 =\sum_i \partial_{h}f_{ijk}(x=h) =\sum_i [(\mu_{ijk}-h)/\sgp_t +\partial_h g_{ijk}]  $
Since $\tau_{4jk}\tauh_{4jk}=\tau_{1jk}\tauh_{1jk}+2 \vt_{jk}$ and  $\tau_{3jk}\tauh_{3jk}=\tau_{2jk}\tauh_{2jk}+2 v_{jk}$, we have 
\BEQ \label{sumAz}
\sum_{i=1}^4 s_i A_{ijk} \partial_x f_{ijk}(h) = \frac{f_{1jk}(h)}{\sgg_t} \sum_{i=1}^4 s_i A_{ijk} (\mu_{ijk}-h) 
\NEQ
Simplifying
\begin{align} \label{sumAzSimp}
\frac{1}{2}\sum_{i=1}^4 s_i A_{ijk} (\mu_{ijk}-h) &= \left[\vt_{jk}(\mu_{4jk}-h)-v_{jk}(\mu_{3jk}-h)\right]  =\nonumber \\
\left[\vt_{jk}(h-ct-\vt_{j,k}\Delta)-v_{jk}(h+ct-2ht-v_{jk}\Delta)\right]  
&=\left[(v_{jk}^2-\vt_{jk}^2)\Delta- ct (v+\vt_{jk}) +h(\vt_{jk}-v_{jk}) + 2ht  v_{jk}\right]  \nonumber\\
=\left[16jkt(1-t)\Delta - 4cjt(1-t)- 4hkt +2htv_{jk}\right]
&=\left[16jk\Delta - 4cj +4h(j-k) \right] t(1-t)
\end{align}
This precisely cancels with $\Bb_{jk}=\sum_i s_i B_{ijk}$ from \eqref{BbCalc}.
For the lower boundary,
we regroup the sum using $j\rightarrow \jh=j+1$, $k\rightarrow \kh=k-1$. This corresponds to centering the generator relative to $\ell$ instead of $h$.
\qed




\subsection{Proof of Corollary \ref{corM0} }\label{AppM0}

\Prf For $M_0$, the coefficients in \eqref{HLCM1g} satisfy $a_{ijk}^{(0)}= \ah_{ijk}^{(0)}= 0$ and $e_{ijk}^{(0)}=1 \Gam_{ijk}$.
Only the last term in \eqref{HLCM1g} is nonzero and the sum reduces to
\begin{align*} \label{CHLMomM0}
M_0 &= \sum_{ijk}\frac{s_{i}}{\sqrt{2\pi}\sig} \Gamma_{ijk}e^{-g_{ijk}}[E_{\sigma_t}(h-\mu_{ijk})- E_{\sigma_t}(\ell-\mu_{ijk})] \\
&= \sum_{jk}\sum_{i=1}^2 \frac{s_{i}}{\sqrt{2\pi}\sig} \Gamma_{ijk}e^{-g_{ijk}} [E_{\sigma_t}(\ell-\mu_{ijk}+2\Delta)- E_{\sigma_t}(\ell-\mu_{ijk})] \ .
\end{align*}
Here $E_{\sigma_t}$ is the scaled $erf$ function, $E_{\sigma_t}(x)\equiv .5* erf(x/\sqrt{2}\sigma_t)$.
To simplify the first sum,  we used $h-\mu_{4jk}= \mu_{1jk}-h$, $h-\mu_{3jk}= \mu_{2jk}-h$,
$E_{\sigma_t}(\ell-\mu_{4jk}) = -E_{\sigma_t}(\ell-\mu_{1jk} + 2 \Delta)$,
$E_{\sigma_t}(\ell-\mu_{3jk}) = -E_{\sigma_t}(\ell-\mu_{2jk} + 2 \Delta)$. The $\Gamma_{ijk}$ satisfy
$\Gamma_{1jk}=(2g_{1jk}-1)w_{jk}^2/\sgg=w_{jk}^2[(c-w_{j,k}\Delta)^2-\sgp]/\sig^4$ and 
$ \Gamma_{2jk} =(2g_{2jk}-1)\wt_{jk}(\wt_{jk}-2)/\sgg=\wt_{jk}\wt_{j,k+1}[(2h-c-\wt_{j,k}\Delta)^2-\sgp]/\sig^4$.

To sum these terms, we reparametrize $k(j,\kh)$. For $i=1,4$, we set $k=\kh-j$, $\kh=k+j$. For $i=2,3$, we set $k=j-\kh$, $\kh=k-j$. With these transformations,
$v_{j,\kh}= j-\kh t$, $\vt_{j,\kh}= j-\kh t$, $w_{j\kh}= 2 \kh$, $\wt_{j\kh}= 2 \kh$, $\mu_{1,j\kh} =ct+\vt_{j,\kh}\Delta$, 
$g_1= (c-w_{j,\kh}\Delta)^2/2\sgg$, $\mu_2= (2h-c)t+v_{j,\kh}\Delta$, $g_2=(2h-c-\wt_{j,\kh}\Delta)^2/2\sgg$.
Since the $g_i$ depend only on $\kh$ and not $j$, so do the $\Gamma_{.j\kh}$.
The double sum splits into a single sum
\BEQ \label{CHLMomM0a}
\sum_{i\kh}\Gamma_{i\kh}e^{-g_{ij\kh}}\sum_j [E_{\sigma_t}(\ell-\mu_{ij\kh}+2\Delta)- E_{\sigma_t}(\ell-\mu_{ij\kh})]=\sum_{\kh}\Gamma_{1\kh}e^{-g_{1\kh}} - \Gamma_{2\kh}e^{-g_{2\kh}}
\NEQ
where we have dropped the $j$ dependence on $g$ and $\Gamma$. We use that for a given $k$, the sum of the integrals for $(1,j,\kh)$
and $(4,j,\kh)$ cover the region from $-\infty$ to $\infty$. This allows us to collapse the sum over $i \in(1,4),j$.
Similarly, the sums over $(2,j,\kh)$ and $(3,j,\kh)$ collapse. We recognize the expression in \eqref{CHLMomM0a} to precisely correspond to $M_0(t,h,\ell,c)= p(h,\ell,c)$ as given by \eqref{dChoiRoh}.
\qed

We would very much like to have expressions for the first and second moment that reduce the double sum to a single sum. This does not appear possible because the $a_{ijk}$ and $a'_{ijk}$ do not vanish. 



\ACKNO{The author thanks the referee for his comments.}
\end{document}